\documentclass[a4paper,10pt]{amsart}

\usepackage{amsmath, amscd, amsthm, amssymb, mathrsfs}
\usepackage[TS1,T1]{fontenc}

\usepackage[varg]{newtxmath}
\usepackage{tgtermes}
\usepackage[scale=0.87]{tgheros}

\DeclareMathAlphabet{\mathsf}{OT1}{qhv}{m}{n}
\SetMathAlphabet{\mathsf}{bold}{OT1}{qhv}{b}{n}

\usepackage{enumerate}

\DeclareMathAlphabet{\emathcal}{U}{esstixcal}{m}{n}

\usepackage[hyphens]{url}
\usepackage[colorlinks,citecolor=blue,urlcolor=blue,linkcolor=blue,breaklinks,bookmarks]{hyperref}
\usepackage[hyphenbreaks]{breakurl}  %

\usepackage[elide]{natbib}

\theoremstyle{remark}
\newtheorem{Remark}{Remark}[section]

\let\ts=\thinspace

\newcommand\la{\langle}
\newcommand\ra{\rangle}

\newcommand{\idsf}{\textup{\textsf{\textbf=}}}
\newcommand{\st}{\ensuremath{\mathsf{t}}}
\newcommand{\si}{\ensuremath{\mathsf{i}}}
\newcommand{\sa}{\ensuremath{\mathsf{a}}}
\newcommand{\ska}{\textup{\textsf{\.{a}}}}
\newcommand{\se}{\ensuremath{\mathsf{e}}}
\newcommand{\ske}{\textup{\textsf{\.{e}}}}
\newcommand{\skke}{\textup{\textsf{\"{e}}}}
\newcommand{\so}{\ensuremath{\mathsf{o}}}
\newcommand{\ex}{\ensuremath{\mathsf{ex}}}

\newcommand{\sis}{\ensuremath{\varepsilonup}} 
\newcommand{\sist}{\ensuremath{\epsilonup}} 

\newcommand{\Lcal}{\textup{\textbf{Ł}}}
\newcommand{\Sh}{\textup{\textbf{Sh}}}
\newcommand{\Shis}{\textup{\Sh{\sis}}}

\newcommand{\For}{\mathrm{For}}
\newcommand{\ForL}{\For_{\Lcal}}

\DeclareFontFamily{U}{bboldvvii}{}
\DeclareFontShape{U}{bboldvvii}{m}{n}
{<-6> bbold5 <6-8> bbold7
<8-> bbold10}{}

\DeclareMathAlphabet{\mathbbx}{U}{bboldvvii}{m}{n}

\newcommand{\uU}{\ensuremath{\mathbbx{U}}}
\newcommand{\DD}{\ensuremath{\mathbbx{D}}}
\newcommand{\dd}{\ensuremath{\mathbbx{d}}}
\newcommand{\dc}{\ensuremath{\mathbbx{c}}}

\newcommand{\mM}{\ensuremath{\mathfrak{M}}}

\newcommand{\cS}{\ensuremath{\mathcal{S}}}
\newcommand{\cP}{\ensuremath{\mathcal{P}}}
\newcommand{\cM}{\ensuremath{\mathcal{M}}}
\newcommand{\cQ}{\ensuremath{\mathcal{Q}}}
\newcommand{\ea}{\ensuremath{\emathcal{a}}}
\newcommand{\eb}{\ensuremath{\emathcal{b}}}

\MakeRobust{\ref}

\makeatletter
\newcommand{\labeltext}[2]{%
  \@bsphack
  \csname phantomsection\endcsname 
  \def\@currentlabel{#1}{\label{#2}}%
  \@esphack}
\makeatother

\begin{document}

\title{The calculus of names -- the legacy of Jan Łukasiewicz}

\author{Andrzej Pietruszczak}

\address{Department of Logic, Institute of Philosophy, Nicolaus Copernicus University in Toruń}

\curraddr{Stanisława Moniuszki 16/20, 87-100 Toruń, Poland}
\email{pietrusz@umk.pl}
\thanks{The research has been supported by the grant from the National Science Centre (NCN), Poland, project no.\ 2021/43/B/HS1/03187.}



\begin{abstract}
\noindent With his research on Aristotle’s syllogistic, Jan  \citet{Luk34, Luk39, Luk57,Luk63} initiates the branch of logic known as the calculus of names. This field deals with axiomatic systems that analyse various fragments of the logic of names, i.e., that branch of logic that studies various forms of names and functors acting on them, as well as logical relationships between sentences in which these names and functors occur. In this work, we want not only to present the genesis of the calculus of names and its first system created by Łukasiewicz. We also want to deliver systems that extend the first. In this work, we will also show that, from the point of view of modern logic, Łukasiewicz’s approach to syllogistic is not the only possible one. However,  this does not diminish Łukasiewicz’s role in the study of syllogism. We believe that the calculus of names is undoubtedly the legacy of Łukasiewicz.
\end{abstract}

\maketitle

\vspace*{-12pt}

\section*{Introduction}

In this work, we want not only to present the genesis of the calculus of names and the first system developed by Jan \citet{Luk34, Luk39, Luk57,Luk63} but also to present systems that are an extension of that initial one, including those enriched with singular sentences of Stanisław Leśniewski’s Ontology, which are not classified as syllogistic. In this work, we will also show that, from the point of view of modern logic, Łukasiewicz’s approach to syllogisms is not the only possible one. It in no way diminishes Łukasiewicz’s role in the study of syllogistic. We believe that the calculus of names is indisputably Łukasiewicz’s legacy.

In the first section, we will present the logic of names and so-called traditional logic. We will present various possible interpretations and forms of categorical sentences in the modern logic of names.

\allowdisplaybreaks

Section~\ref{sec2} will be devoted to the calculus of names as a specific development of traditional logic. We will present the origins of this calculus and Łukasiewicz’s original system. We give the set-theoretic semantics of this system and show that it is equivalent to the lexical semantics (when we take an appropriate set of non-empty general names for the name variables). We will note that Łukasiewicz's calculus is \emph{sound} and \emph{complete} not only with respect to the set of all non-empty general names but also with respect to the set of general names having at least two references (we will formally prove these facts in Section~\ref{sec6}).

In Section~\ref{sec3}, we present other possible approaches to formalising Aristotle's syllogistic. They come from, among others, John \citet{Cor72a, Cor74}, Timothy John \citet{Smil} and Robin \citet{Sm}. We also present a “competitive approach” in the form of the sequent calculus.

In Section~\ref{sec4}, we will discuss two modern takes on the calculus of names that allow them to be applied to empty names. The first of them (using the so-called weak interpretation of universal affirmative sentences) was reported and researched by J.\thinspace C.~\citet{Sh}. We will present a series of definitional extensions of his system and also give its set-theoretic semantics. We will note that Shepherdson's system is \emph{sound} and \emph{complete} not only with respect to the set of all general names but also with respect to the set of general names that are either empty or have at least two references, i.e., we exclude names having exactly one reference (we will formally prove these facts in Section~\ref{sec6}). The second approach (using the so-called strong interpretation of universal affirmative sentences) was initiated by Jerzy \citet{Sl}. However, his system is not complete. It was noted in \citep{ja87}, where complete extensions of Słupecki’s system are also given.

In Section~\ref{sec5}, we present and analyse extensions of both types of systems with singular sentences of Leśniewski's Ontology. Arata \citet{Ish} gave the propositional (quantifier-free) fragment of this theory. Firstly, we analyse the fusion of Shepherdson's system with this fragment. We give four axiomatisations of this fusion. Moreover, we will present a series of definitional extensions of this fusion and also give its set-theoretic semantics. We will note that it is sound and complete with respect to the set of all general names (see Section~\ref{sec6}). Secondly, we analyse the fusions of complete extensions of Słupecki's system with the quantifier-free fragment of Ontology. We give axiomatisations of these fusions and show that they are definitionally equivalent to the fusion of Shepherdson's system with the quantifier-free fragment of Ontology.

In Section~\ref{sec6}, we present different approaches to the proof of completeness of calculi of names with respect to set-theoretic semantics. The first approach comes from \citep{Sh}, where is used a technique similar to that used to prove Stone's representation theorem for the elementary theory of Boolean algebras. This approach uses appropriate filters constructed from the elements of a given algebra (a model of a given theory). The second approach consists of the appropriate direct application of Henkin's method to calculi of names. This method is commonly used in the proof of the completeness of the propositional logic or predicate logic. In it, we use canonical models built for maximal consistent sets in a given calculus. We give two ways of doing this.

In the last section, we briefly present other possible extensions of the systems considered earlier by adding a few new kinds of sentences. They will be traditional singular sentences, Czeżowski’s singular sentences (with a subject of the form ‘this~$S$’) and identities for singular names.

\section{The logic of names and traditional logic\label{sec1}}

\subsection{The logic of names}

The logic of names is constructed using the method of logical schemes  \citep[see, e.g.,][]{Q50}, which consists of the fact that, based on the analysis of the surface syntactic structure of sentences and expressions of natural language, sentence schemes are introduced in which various types of schematic letters appear instead of names. Name logic is an intermediate link between propositional logic and predicate logic. In propositional logic, we study relations between sentences but are not interested in the syntactic structure of sentences in which there are no propositional conjunctions. We look only at those relationships that depend solely on propositional connectives. In quantifier logic, the opposite is true; we analyse the deep structure of sentences using quantifiers binding variables and additionally introduced sentential connectives. These are characteristics of the modern mathematical stage of formal logic.

The logic of names can be regarded as a systematic development of specific fragments of traditional, pre-mathematical formal logic. We include not only the known piece of it, which is Aristotle's syllogistic, but also studies on compound names and relative names. The latter, as oblique syllogisms, were already considered by Aristotle in his Prior Analytics and Joachim Jungius in his Logica Hamburgensis. However, a systematic theory of them appeared only in the nineteenth century in the works of Hamilton, Schröder and de Morgan. De Morgan analysed reasoning such as: since every horse is a mammal, then every horse's head is a mammal's head.

Aristotelian syllogistic deals with categorical sentences having one of the following four schemes:
\begin{enumerate}[\textbullet]
\item Every $S$ is a $P$\hfill	(\emph{universal affirmative})
\item Some $S$ is a $P$ \hfill	(\emph{particular affirmative})
\item No $S$ is a $P$	\hfill	(\emph{universal denial})
\item Some $S$ is not a $P$ \hfill	(\emph{particular denial})
\end{enumerate}
These schemes can be used for any general name. The letter ‘$S$’ is to be replaced by a general name appearing as a subject, and the letter ‘$P$’ is to be replaced by a general name appearing as a predicate. The following natural interpretation of these sentences is generally accepted:
\begin{enumerate}[\textbullet]
\item ‘Every $S$ is a $P$’ is true if and only if the extension of the name $S$ is included in the extension of the name $P$;
\item `Some $S$ is a $P$' is true if and only if the names $S$ and $P$ have a common referent;
\item `No $S$ is a $P$' is true if and only if the extensions of the names $S$ and $P$ are disjoint;
\item `Some $S$ is not a $P$' is true if and only if the name $S$ has a referent that is not a referent of the name~$P$.
\end{enumerate}
To avoid any further ambiguity, we emphasise that the extension of a given general name is the distributive set of all its referents (i.e., it is its extension in Frege's sense). Thus, every empty general name (i.e., this without any referent) has as its extension the empty set (which is included in every set).

The method of logical schemes allows the study of a broad class of natural language sentences. In addition to categorical sentences, we can also study singular sentences of the form ‘$a$ is a $P$’ and ‘$a$ is not a $P$’, in the subject of which there is a name pointing to exactly one object. We also have a whole spectrum of sentences corresponding to categorical sentences; for example, these are sentences like ‘$S$ is the same as a $P$’ (or otherwise: ‘Every $S$ is a $P$ and vice versa’, ‘All $S$ is a $P$ and vice versa’), ‘Exactly one $S$ is a $P$’, ‘At most one $S$ is a $P$’, ‘The only $S$ is a $P$’, and others. We can analyse plural sentences such as ‘Exactly two $S$s are $P$s’, ‘At least two $S$s are $P$s’, ‘At most two $S$s are $P$s’, and others. We can also treat modal versions of these sentences in which the copula ‘is’ is replaced by one of the phrases: ‘must be’, `maybe’; the phrase `is not’ is replaced with one of the phrases `must not be’, `must not be’, ‘may not be’, ‘cannot be’. It is also possible to analyse sentences whose subjects and predicates have compound names of the form: ‘$S$ and $P$’, ‘$S$ or $P$’, and ‘not-$S$’. The same applies to relative names such as ‘friend’, ‘mother’, and others. We can also transform the latter sentences from active to passive (e.g., ‘reader’ to ‘read by’) and, from two such names, create a third relative name (e.g., ‘mother’s father’). We do not have to limit ourselves to relative terms but may extend the approach to verbs, e.g., instead of `is a reader’, we take ‘reads’ \citep[see, e.g.,][]{M,P-HiM}.

Rich metalogical research can be carried out on the material given above: various types of set-theoretic semantics and axiomatisations of different fragments of the logic of names consistent with it can be introduced.

\subsection{Traditional logic\label{subsec:TL}}

For some reasons \citep[{cf.}][pp.~103--104]{Cor74} Aristotelian syllogistic was applied only to non-empty general names (i.e., those with at least one referent). Traditional logic was a continuation of Aristotle's syllogistic. Therefore, it also took over the limitation to apply only to non-empty general names. Traditional logic, by argument schemes, studied the logical relationships between categorical sentences. The primary connection between them is the entailment relation (in symbols: $\therefore$). It is the converse of logical consequence. In traditional logic, the entailment relation was expressed by so-called correct argument schemas (or: argument forms). If categorical propositions are understood naturally, when we limit the applications to non-empty names, all the argument schemes distinguished by traditional logic are valid because true premises always give a true conclusion. So, in traditional logic, the argument schemas below are valid:
\begin{enumerate}[\textbullet]
\item Every $S$ is a $P$  $\therefore$ Some $S$ is a $P$ 			                \hfill \emph{subalternation}
\item Every $S$ is a $P$  $\therefore$ Some $P$ is an $S$ \hfill		               \emph{conversion per accidens}
\item No $S$ is a $P$  $\therefore$ Some $S$ is not a $P$ \hfill			              \emph{subalternation}
\item No $S$ is a $P$ $\therefore$ Some $P$ is not an $S$ \hfill		               \emph{conversion per accidens}
\item Every $S$ is a $P$ $\therefore$ It is not the case that no $S$ is a $P$	\hfill  \emph{contrariety}
\item 	It is not the case that some $S$ is a $P$ $\therefore$  Some S is a $P$ \hfill	                \emph{subcontrariety}
\item 	Every $M$ is a $P$, Every $S$ is a $M$ $\therefore$  Some S is a P	\hfill		            \emph{Barbari}
\item 	Every $M$ is a $P$, Every $M$ is a $S$ $\therefore$   Some S is a $P$ \hfill		 	             \emph{Darapti}
\item 	Every $P$ is a $M$, Every $M$ is a $S$ $\therefore$  Some $S$ is a $P$ \hfill			            \emph{Bamalip}
\item 	No $M$ is a $P$, Every $S$ is a $M$ $\therefore$  Some S is not a $P$ \hfill		           \emph{Celaront}
\item 	No $P$ is a $M$, Every $S$ is a $M$ $\therefore$ Some S is not a $P$ \hfill		              \emph{Cesaro}
\item 	Every $P$ is a $M$, No $S$ is a $M$ $\therefore$  Some S is not a $P$ \hfill		         \emph{Camestros}
\item 	Every $P$ is a $M$, No $M$ is a $S$ $\therefore$  Some S is not a $P$ \hfill		            \emph{Calemos}
\item 	No $M$ is a $P$, Every $M$ is a $S$ $\therefore$ Some $S$ is not a $P$ \hfill		           \emph{Felapton}
\item 	No $P$ is a $M$, Every $M$ is a $S$ $\therefore$ Some $S$ is not a $P$ \hfill		              \emph{Fesapo}
\end{enumerate}
Moreover, the scheme ‘Some $S$ is an $S$’ is generally true (i.e. it is valid). Notice that if we allow empty general names, the above argument schemas and the sentence schema are no longer valid.

\subsection{A contemporary approach to the logic of names\label{subsec:1.3}}

In the contemporary approach, we allow empty general names, so some argument forms of traditional logic are no longer valid. So, the following question arises:
\begin{enumerate}[\textbullet]
\item Can the meaning of categorical sentences be changed to preserve the validity of argument forms of traditional logic, even when substituting empty names is allowed?
\end{enumerate}
This new interpretation is to meet, however, the following condition:
\begin{enumerate}[\textbullet]
\item When terms are limited to non-empty terms, it coincides with natural usage.
\end{enumerate}

Therefore, to “save” the subalternation and conversion per accidens of affirmative sentences, contrariety and syllogisms \emph{Barbari}, \emph{Darapti}, \emph{Bamalip}, \emph{Celaront}, \emph{Cesaro}, \emph{Camestros}, \emph{Felapton} and \emph{Fesapo}, the so-called \emph{strong} interpretation of universal affirmative sentences was used, where for its truth, we require the non-emptiness of the name in its subject. Thus:
\begin{enumerate}[\textbullet]
\item A universal affirmative sentence (in the strong interpretation) is true if and only if it has a non-empty name in the subject, which extension is included in the extension of the name from the predicate.
\end{enumerate}
An interpretation where we do not apply the added requirement is called \emph{weak}. Of course, for non-empty names, the two interpretations are indistinguishable.

Note that with the strong interpretation, if we allow empty names, then:
\begin{enumerate}[\textbullet]
\item 	‘Some $S$ is not a $P$’ is not a contradiction of ‘Every $S$ is a $P$’ and vice versa.
\end{enumerate}
Indeed, if $S$ is an empty name, both sentences are false.

For this reason, Tadeusz \citet[pp.~233--234]{Kot}, followed by Czesław \citep[pp.~128--130]{Lej}, proposed the use of additional universal affirmative sentences of the form ‘All $S$ is a $P$’, which are true for any empty name standing in the subject (regardless of what name stands in the predicate). These sentences, therefore, have the interpretation of universal affirmative sentences in the weak sense. Kotarbiński and Lejewski’s proposal indicates that they believed that in the meaning of the phrase ‘all~$S$’, there is no reservation about the non-emptiness of $S$; that is, this reservation is implicitly related to the phrase ‘every~$S$’. When we limit ourselves to non-empty names, the interpretations of both types of universal affirmative sentences coincide. Namely, what is supposed to be implicit in the meaning of ‘every~$S$’ is implied in the assumption imposed on the names.

Note that for new universal affirmative sentences, for all general names:
\begin{enumerate}[\textbullet]
\item 	‘Some $S$ is not a $P$’ is a contradiction of ‘All $S$ is a $P$’ and vice versa.
\end{enumerate}

To “save” the subalternation of denial sentences and syllogism \emph{Camestros}, it is enough to use the \emph{strong} interpretation for universal denial sentences, where we require the non-emptiness of the name standing in the subject for their truth. Thus:
\begin{enumerate}[\textbullet]
\item 	A universal denial sentence (in the strong interpretation) is true if and only if it has a non-empty name in the subject whose extension is disjoint with the extension of the name of the predicate.
\end{enumerate}
\citet[p.~130]{Lej} proposed the introduction of two functors to construct universal denial sentences. In addition to the functor of weak exclusion ‘no \dots\ is \dots’, he introduced the functor of strong exclusion ‘every \dots\ is not \dots’.

To “save” the \emph{conversion per accidens} of denial sentences and syllogism \emph{Calemos}, we must use an even stronger, \emph{super-strong} interpretation of universal denial sentences, requiring both names to be non-empty for their truth. Thus:
\begin{enumerate}[\textbullet]
\item A universal denial sentence (in the super-strong interpretation) is true if and only if it has non-empty names in both the subject and the predicate, the extensions of which are disjoint.
\end{enumerate}
Neither Kotarbiński nor Lejewski used this interpretation. Furthermore, they did not introduce new sentences expressing it. In \citep{ja87,ja91b, ja91c}, it was proposed that these should be sentences of the form ‘Every $S$ is not a $P$ and vice versa’. It was modelled on Kotarbiński’s comments on the phrase ‘every~$S$’ and on the sentences he used in the form ‘All $S$ is a $P$ and vice versa’ (which state that the ranges of both names are equal).

Again, when we limit ourselves to non-empty general names, the interpretations of the three types of universal denial sentences coincide. Indeed, this is implicit in the meaning of ‘every~$S$’ and explicitly implied in the assumption imposed on the names. Moreover, what is expressly contained in the meaning of the phrase ‘and vice versa’ is implicit in interpreting the functor ‘no \dots is \dots’.

\section{Calculus of names as an extension of traditional logic\label{sec2}}
\subsection{The genesis of the calculus of names\label{subsec:2.1}}

As stated in the introduction, Łukasiewicz is undoubtedly the creator of the calculus of names. The following words from him \citep{Luk34} show the genesis of this calculus:\footnote{The Polish text of \citep{Luk34} is translated by the author of this paper. The Polish term ‘reguła wnioskowania’ is translated as ‘rule of inference’.}
\begin{quotation}\small
\textbf{5.} A fundamental difference exists between a logical thesis and a rule of inference.

A \emph{logical thesis} is a sentence in which, apart from logical constants, there are only sentence or name variables, which is true for all values of the variables that occur in it. An inference rule is a prescription that authorises a person which make inferences to derive new theses based on recognised theses. For example, [\dots] the principles of identity\label{identity} [such as “If $p$, then $p$” and “Every $a$ is an $a$”] are logical theses, but the rule of inference is the following «rule of detachment»:

Whoever accepts as true the implication “If $\alpha$, then $\beta$” and the antecedent of this implication “$\alpha$” has the right to accept as true also the consequent of this implication~“$\beta$”.
\end{quotation}
The problem is that the fact that a given implication is considered to be true can be understood differently.  Since the antecedent $\alpha$ and the consequent $\beta$ appearing in the implication considered by Łukasiewicz have variables\footnote{Here, ‘variable’ may be replaced by ‘schematic letter’. Willard Van Orman \citet{Q53} wrote about the significant difference between these terms. The name ‘calculus of names’ probably comes from the fact that it is about “variables” for which we substitute names and perform some “calculus” on its formulas.}, two situations can be considered:
\begin{enumerate}[1.]
\item If with a given admissible substitution for variables (schematic letters), the schemes ‘If $\alpha$, then $\beta$’ and $\alpha$ give true sentences, then, with this substitution, we also have a true sentence obtained from the schema $\beta$.

\item If the formulas ‘If $\alpha$, then $\beta$’ and $\alpha$ are logically valid (i.e., they give true sentences under any admissible substitution for variables), then $\beta$ is also logically valid.
\end{enumerate}
Thus, in the first of the above points, we treat the rule of detachment as the following valid \emph{argument form} (or \emph{argument schema}):
\[
\frac{\text{$p$\qquad\quad If $p$ then $q$}}{q} \qquad\qquad\frac{\;p\qquad\quad p\to q\;}{q}
\]
where the letters ‘$p$’ and ‘$q$’ stand in the place of sentences. However, the second point says that the rule of detachment from two logically valid formulas always leads to such a formula, and so it has the following scheme:
\[
\frac{\text{$\alpha$\qquad\quad If $\alpha$ then $\beta$}}{\beta} \qquad\qquad\frac{\;\alpha\qquad\quad\alpha\to\beta\;}{\beta}
\]
where $\alpha$ and $\beta$ represent arbitrary sentence formulas. Here, this rule is then something that, from two logically valid formulas “produces” a third one. Depending on needs, the rule of detachment can perform one of the above roles or both.

The detachment rule is used in the latter role in logical calculi, including the calculus of names. It is a “generator” of theses, which has to be logically valid. This generator derives theses from axioms. The detachment rule says that if $\alpha$ and $\alpha\to\beta$ are already justified, then $\beta$ is justified. This means that a line in a derivation containing $\beta$ is justified, provided that (for some $\alpha$) both $\alpha$ and $\alpha\to\beta$ appear in the derivation before~$\beta$. Since initial theses (axioms) are logically valid, and the rules used (including the rule of detachment) transform logically valid formulas into new such formulas. Therefore, all these obtained with their help are also logically valid.

However, in some logical calculi, not all rules perform both of the roles indicated above. The primary example is the substitution rule used by Łukasiewicz. It says that from any logically valid formula, every permissible substitution for “variables” gives a logically valid formula. This rule does not even have a scheme by which to express it. Hence, it cannot be “confused” with a scheme of correct reasoning (argument form).

Another example, which has a scheme, is the following \emph{rule of necessitation}:
\[
\frac{\alpha}{\text{\;It is necessary that~}\alpha\;}\qquad\qquad \frac{\;\alpha}{\Box\alpha\;}
\]
It is accepted in all normal modal logics. It takes us from a logically valid formula to a new logically valid formula. However, the following argument form is not valid (where the letter ‘$p$’ stands in the place of a sentence):
\[
\frac{p}{\;\text{It is necessary that~}p\;}\qquad\qquad\frac{p}{\;\Box p\;}
\]
Indeed, we have true sentences that are not necessary. For this reason, we cannot reason according to this scheme in modal theories, where we are interested in true sentences which are not logically true.

The above remarks generally refer to the genesis of various types of logical calculi. The origin of the calculus of names itself can be found by continuing the quote from \citep{Luk34}:
\begin{quotation}\small
\textbf{6.} The original Aristotelian syllogism is a logical thesis, the traditional syllogism has the meaning of a rule of inference.

The Barbara mode given [below], [\dots], is an implication of the type “If $\alpha$ and $\beta$, then $\gamma$”, [\dots]. As an implication, an Aristotelian syllogism is a proposition that Aristotle holds to be true, namely that the proposition is true for all values of the variables “$a$”, “$b$” and “$c$” that occur in it. Therefore, we get true sentences if we substitute some constant values for these variables. Since in the considered mode, apart from variables, there are only logical constants, namely “if-then”, “and” and “every-is”, the Aristotelian syllogism is a logical thesis.

The traditional syllogism:

$\dfrac{\parbox{20mm}{Every $b$ is an $a$\\ Every $c$ is a $b$}}{\text{Every $c$ is an $a$}}$

\noindent is \emph{not} an implication. It consists of three sentence forms, listed one under the other, which do not form a single sentence. Since a traditional syllogism is not a proposition, it cannot be true or false either since, according to the generally accepted view, truth and falsity belong only to propositions. A traditional syllogism is, therefore, \emph{not} a thesis. If we substitute some constant values for the variables in this syllogism, we do not get a proposition but an \emph{argument form}. So a traditional syllogism is an argument schema and has the meaning of a \emph{rule of inference}, which can be more precisely expressed as follows:

Whoever accepts as true premises of the form ``Every $b$ is an $a$'' and ``Every $c$ is a $b$'' has the right to accept as true a conclusion of the form ``Every $c$ is an $a$''.\footnote{Footnote 11 added: ``How imprecise the historical studies of logic to date are is evidenced by this very characteristic detail: all the authors I know who have written about Aristotelian logic, [\dots], present Aristotelian syllogisms in the traditional form, without even realising the fundamental difference between these forms.''}

\textbf{7.} Thanks to the distinction between logical theses and rules of inference, it became possible for logical sciences to construct axiomatically in the form of deductive systems.
\end{quotation}
The problem, however, is that Łukasiewicz used the term ‘rule of inference’ in two meanings. Firstly, he writes that a traditional syllogism is a “rule of inference” qua a valid argument form (or a valid argument schema), i.e., it preserves the truth in the following sense: if its premises are true, then its conclusion is also true. Secondly, the detachment rule he used in his calculus of names is a “rule of inference” preserving validity, i.e., from two logically valid formulas, always leads to such a formula. It is a “generator” of theses that must be logically valid since used axioms are logically valid.

Łukasiewicz states that syllogisms cannot be treated as rules of inference because such rules can transform only sentences. However, rules of inference can be viewed differently. Namely, as in the sequent calculus, they can be viewed as transforming argument forms (see Section~\ref{sec3}). It does not mean we think that Aristotle’s syllogistic should be treated as a sequent calculus, though the system presented in the next subsection can be transformed into a sequent calculus (see Section~\ref{sec3}). However, we mainly wanted to show that the term ‘rule of inference’ can be used in a third sense.

\subsection{Łukasiewicz's calculus of names\label{subsec:2.3}}
Łukasiewicz presented his reconstruction of Aristotle's syllogistic as a calculus of names for the first time in \citep{Luk63} (the first Polish edition in 1929). He repeated it in his works \citep{Luk34,Luk39} and then in \citep[{\S}25; the first edition in 1951]{Luk57}. This reconstruction is presented in the continuation of the previously quoted text from \citep{Luk34}:
\begin{quotation}\small
\textbf{8.} The theory of the Aristotelian syllogism, which Aristotle has already tried to axiomatise, but which has not yet been presented in an axiomatic form, is based on two fundamental concepts: “Every $a$ is a $b$”, in the signs “$Uab$”, and “Some $a$ is a $b$”, in the signs “$Iab$” and on the following axioms:

1. Every $a$ is an $a$.

2. Some $a$ jest an $a$.

3. If every $b$ is an $a$ and every $c$ is a $b$, then every $c$ is an $a$.

4. If every $b$ is an $a$ and some $b$ is a $c$, then some $c$ is an $a$.

In the signs (the functors ``$U$'' and ``$I$'' come before the arguments, and such same the conjunction sign ``$K$'' = ``and''):

1. $Uaa$.

2. $Iaa$.

3. $CKUbaUcbUca$ (\emph{Barbara}).

4. $CKUbaIbcIca$ (\emph{Datisi}).

The expressions ``Some $a$ is not a $b$'', in the signs ``$Oab$'', and ``No $a$ is a $b$'', in the signs ``$Yab$''\footnote{Łukasiewicz used his bracketless notation here. In \citep{Luk57} `$U$' and `$Y$' were replaced by `$A$' and `$E$', respectively.}, can be defined as follows:

Df1. $Oac = NUab$.

Df2. $Yab = NIab$.

By both rules of substitution and detachment (propositional variables may be substituted with propositional forms of Aristotelian logic, for name variables \emph{only} other name variables), and with the help of theses of propositional logic, from these axioms and definitions, we can derive all 24 (not 14 nor 19!) the correct modes of Aristotelian syllogistic.\footnote{Footnote 14 added: “The axiomatisation of Aristotelian syllogistic presented here, as well as the deduction of all modes, can be found in the script from my lectures, delivered in the autumn trimester of 1928/29 at the University of Warsaw, entitled: \emph{Elementy logiki matematycznej} [\dots]” \citep[{cf.}][]{Luk63}).}
\end{quotation}

We will use the \textsf{a}a following abbreviations for schemes of categorical sentences:\footnote{The abbreviations are derived from the vowels in the Latin words `\emph{affirmo}' and `\emph{nego}'.}
\begin{tabbing}
\mbox{\qquad}\= $S\sa P$ \= --X \= for `Every $S$ is a $P$' \kill
\> $S\sa P$ \> ~-- \> for `Every $S$ is a $P$'\\
\> $S\si P$ \> ~-- \> for `Some $S$ is a $P$'\\
\> $S\se P$ \> ~-- \> for `No $S$ is a $P$'\\
\> $S\so P$ \> ~-- \> for `Some $S$ is not a $P$'
\end{tabbing}
Using the above abbreviations, we will reconstruct the Łukasiewicz calculus, which we will denote by \Lcal. We will use the countably infinite set GN of name letters (for which we use ‘$S$’, ‘$P$’, ‘$M$’ and ‘$Q$’ with or without indices). From the name letters and the constants ‘$\sa$’, ‘$\si$’, ‘$\se$’ and `$\so$’, we build \emph{atomic formulas} corresponding to the abbreviations of categorical sentences given above. Moreover, we use the Boolean propositional connectives ‘$\neg$’, ‘$\wedge$’, ‘$\vee$’, ‘$\to$’ and ‘$\leftrightarrow$’ (for negation, conjunction, disjunction, material implication and equivalence) and brackets. The set of atomic formulas determines the set $\ForL$ of all formulas of \Lcal, which is built in the standard way from the atomic formulas, the Boolean propositional connectives and brackets. Thus, $\ForL$ is the smallest set $X$ such that:
\begin{enumerate}[\textbullet]
\item all atomic formulas belong to $X$,
\item if $\alpha,\beta\in X$ and ${\ast}\in\{\wedge, \vee, \to,  \leftrightarrow\}$, then $\neg\alpha\in X$ and $(\alpha\ast\beta)\in X$.
\end{enumerate}

Let us adopt the following notation of Łukasiewicz’s four axioms. The first two are the principles of “identity” \citetext{see the quote from \citealp{Luk34}, on p.~\pageref{identity}}, and the next two correspond to traditional syllogisms \emph{Barbara} and \emph{Datisi}:
\begin{gather*}
S\sa S \label{Ia}\tag{I\sa} \\
S\si S \label{Ii}\tag{I\si} \\
(M\sa P\wedge S\sa M)\rightarrow S\sa P \label{Barbara}\tag{Barbara} \\
(M\sa P\wedge M\si S)\rightarrow S\si P \label{Datisi}\tag{Datisi}	
\end{gather*}
In the last quote, Łukasiewicz adopted three rules for deriving theses: detaching, substituting and defining. The last one does not apply today. Instead, additional specific axioms are introduced as equivalences called \emph{definitions}.\footnote{It can also be assumed that the notations of the form `$S\se P$' and `$S\so P$' are only abbreviations for the formulas `$\neg\, S\si P$' and `$\neg\, S\sa P$', respectively. It means that the formers are not among the formulas at all. Such a solution for `$S\so P$' was attempted by \citet{Sl}; see further point~\ref{subsec4.2}.} In our case, these equivalences will assume the contradictions of the pairs $S\se P$\thinspace--\thinspace$S\si P$ and $S\sa P$\thinspace--\thinspace$S\so P$:
\begin{gather*}
S\se P\leftrightarrow\neg\,S\si P\label{dfe}\tag{df\,{\se}}\\
S\so P\leftrightarrow\neg\,S\sa P\label{dfo}\tag{df\,{\so}}
\end{gather*}
\indent Łukasiewicz added fourteen tautologies of implication-negation fragment of classical propositional logic (CPL), in which he replaces propositional variables with schemes of categorical sentences and their conjunctions \citep[{cf.}][pp.~88--89]{Luk57}. To facilitate the derivation of theses, all substitutions of all CPL tautologies with formulas of the calculus of names can adopted as axioms. However, the essence of this calculus is contained in its specific axioms, i.e., in \eqref{Ia}, \eqref{Ii}, \eqref{Barbara}, \eqref{Datisi}, \eqref{dfe} and \eqref{dfo}.

Further, to obtain theses of \Lcal, without detailed explanations, we will use the necessary tautologies of CPL and the rules of detachment and uniform substitution. We will only specify on the left side of a given thesis which axioms or previously obtained theses should be used. We obtain all implications corresponding to the argument schemes of traditional logic given in Subsection~\ref{subsec:TL}:

\begin{tabbing}
\= (\emph{Camestros})\qquad \= $(M\sa P\wedge M\sa S)\to S\si P$\qquad \=  \kill
\> $(\mathrm{C}{\si})$ \labeltext{C$\si$}{Ci} \> $P\si S\rightarrow S\si P$ \> \emph{conversion}, by \eqref{Ia} and \eqref{Datisi}\\
\> $(\mathrm{C}{\se})$ \labeltext{C$\se$}{Ce} \> $P\se S\rightarrow S\se P$ \> \emph{conversion}, by \eqref{Ci} and \eqref{dfe}\\
\> $({\sa}\mathrm{S}{\si})$ \labeltext{{\sa}S{\si}}{aSi} \> $S\sa P \rightarrow S\si P$ \>  \emph{subalternation}, by  \eqref{Ii} and \eqref{Datisi}\\
\> \> $S\sa P\to P\si S$ \> \emph{conversion per accidens}, by \eqref{aSi} and \eqref{Ci}\\
\> $({\se}\mathrm{S}{\so})$  \labeltext{{\se}S{\so}}{eSo} \> $S\se P \rightarrow S\so P$ \>  \emph{subalternation}, by  \eqref{aSi}, \eqref{dfe} and  \eqref{dfo}\\
\> \> $S\se P\to P\so S$ \> \emph{conversion per accidens}, by \eqref{eSo} and \eqref{Ce}\\
\> \> $S\sa P\to \neg\, S\se P$ \>  \emph{contrariety}, by \eqref{aSi} and \eqref{dfe}\\
\>\> $\neg\, S\si P\to S\so P$  \> \emph{subcontrariety}, by \eqref{aSi} and \eqref{dfo}\\
\> {\normalfont(\emph{Barbari})} \> $(M\sa P\wedge S\sa M)\to S\si P$  \> by \eqref{Barbara} and \eqref{aSi} \\
\> {\normalfont(\emph{Darapti})} \labeltext{Darapti}{Darapti} \> $(M\sa P\wedge M\sa S)\to S\si P$ \> by \eqref{Datisi} and \eqref{aSi} \\
\> (\emph{Bamalip})  \> $(P\sa M\wedge M\sa S)\to S\si P$ \> by \eqref{Barbara}, \eqref{aSi} and \eqref{Ci} \\
\> (\emph{Camestros}) \> $(P\sa M\wedge S\se M)\to S\si P$ \> by (\emph{Barbari}), \eqref{dfe} and \eqref{dfo}\\
\> (\emph{Calemos}) \> $(P\sa M\wedge M\se S)\to S\si P$ \> by (\emph{Camestros}) and \eqref{Ce}\\
\> \emph{Felapton}) \> $(P\se M\wedge M\sa S)\to S\so P$ \> by (\emph{Barbari}), \eqref{dfe} and \eqref{dfo}\\
\> (\emph{Fesapo}) \> $(M\se P\wedge M\sa S)\to S\so P$ \> by (\emph{Felapton}) and \eqref{Ce}\\
\> (\emph{Celaront}) \> $(M\se P\wedge S\sa M)\to S\so P$  \> by \eqref{Darapti}, \eqref{dfe} and \eqref{dfo} \\
\> (\emph{Cesaro}) \> $(P\se M\wedge S\sa M)\to S\so P$ \> by (\emph{Celaront}) and  \eqref{Ce}\\[-5pt]
\end{tabbing}
Note that apart from \eqref{Ci} and \eqref{Ce}, the rest of the above implications cannot be obtained without using \eqref{Ii}. Moreover, similarly, but without using \eqref{Ii}, we obtain implications corresponding to the remaining thirteen traditional syllogisms.

We consider Łukasiewicz’s system to be traditional because, just like traditional logic, we can apply it only to non-empty names. We say that a formula is a \emph{traditional lexical tautology} if and only if it gives true sentences in all substitutions of non-empty general names for name letters (variables). It is proved \citep[cf., e.g.,][]{ja91c} that this system is \emph{sound} and \emph{complete} in the sense that all lexical tautologies and only them are its theses.

The great advantage of the approach to the logic of names proposed by Łukasiewicz is that for the semantic study of a given system, we can use methods known from the metatheory of propositional and predicate logics \citep[cf., e.g.,][]{ja91a, ja91c, ja92}. Moreover, we can consider a given system one of the open first-order theories and use meta-theorems about such theories \citep[see][]{Sh, ja92}.
\begin{Remark}\label{iLcal}
We can also consider another version of \Lcal, which we obtain by rejecting the rule of uniform substitution and accepting all substitution instances of the axioms of \Lcal\ as its specific axioms. We can show that both versions have the same theses \citep[{cf.}][]{Q50}. This second version was used, among others, in \citep{Sh}. \qed
\end{Remark}

\subsection{Set-theoretic semantics for Łukasiewicz’s calculus\label{subsec2.3}}

\emph{Traditional models.} In the semantic study of \Lcal, from a formal point of view, instead of talking of substitutions of general names for name letters, for $\ForL$, it is better to use set-theoretic semantics, which uses \emph{traditional models} of the form $\la \uU ,\DD \ra$, where:
\begin{enumerate}[\textbullet]
\item \uU\ is a non-empty set (\emph{universe}),
\item \DD\ is a function of \emph{denotation}, which assigns to any name letter a non-empty subset of \uU.\footnote{\label{footD}Like \citet[p.~103]{Cor74}, we could assume that the model (interpretation) is just the mapping $\DD$ itself, which assigns to any name letter a non-empty set.}
\end{enumerate}
Using the previously given interpretation of the functors ‘\sa’, ‘\si’, ‘\se’ and ‘\so’, we introduce the notions of \emph{being a true formula} in a model $\mM=\la \uU ,\DD \ra$. For atomic formulas, for any letters \cS\ and \cP\ we assume:
\begin{enumerate}[\textbullet]
\item $\cS\sa\cP$ is true in \mM\ iff the set $\DD(\cS)$ is included in the set $\DD(\cP)$;

\item$\cS\si\cP$ is true in \mM\ iff the sets $\DD(\cS)$ and $\DD(\cP)$ have a common element;

\item $\cS\se\cP$ is true in \mM\ iff the sets $\DD(\cS)$ and $\DD(\cP)$ have no element in common;

\item $\cS\so\cP$ is true in \mM\ iff $\DD(\cS)$ has an element which is not an element of $\DD(\cP)$.
\end{enumerate}
For formulas built with propositional connectives, we use truth tables, i.e., we interpret these connectives as in CPL.

We say that a formula is a \emph{traditional} (\emph{set-theoretic}) \emph{tautology} (or is \emph{traditionally valid}) if and only if it is true in all traditional models. As can be easily seen, \Lcal\ is \emph{sound} with regard to set-theoretic semantics. Indeed, all axioms of \Lcal\ are traditional tautologies, uniform substitution and the detachment rule preserves the traditional validity of formulas. In Section~\ref{sec6}, we show different ways of proving the completeness of~\Lcal.

\medskip

\paragraph{\emph{Set-theoretic tautologies vs.\ lexical tautologies.}} Since to each general name, we can assign a set that is its extension, every set-theoretic tautology is also a lexical tautology. On the other hand, if every set was an extension of some name, then without any additional conditions, it could be proved that the opposite holds, i.e., that every lexical tautology is a set-theoretic tautology.

In fact, we lack names; we cannot assign a name to each set that would cover all its elements. However, as we know, in the model theory of predicate logic, and therefore logic of names, not all sets are needed. Those that are determined by the formulas of elementary number theory with one free variable will suffice. This theory covers what can be said about natural numbers using the names of individual numbers, addition, multiplication, equal sign, propositional connectives and quantifiers. Both definitions of tautologies denote the same set of formulas if the natural language, whose general names we substitute for name letters, satisfies the following condition:
\begin{enumerate}[\textbullet]
\item For every set of natural numbers described above, there is a general name whose extension is this set.
\end{enumerate}
The above condition is not too high. After all, we are talking about such general names as ‘smallest natural number’, ‘largest natural number’, ‘even number’, ‘number greater than 10’ and similar others.

\medskip

\paragraph{\emph{Polyreferential set-theoretic semantics for Łukasiewicz’s calculus.}} A~general name which has at least two referents (resp.\ exactly one referent) we will call \emph{polyreferential} (resp.\ \emph{monoreferential}). To those names correspond to \emph{polyreferential models} of the form $\la \uU ,\DD \ra$, where \uU\ is a set having at least two elements and \DD\ is a function of \emph{denotation}, which assigns to any name letter a subset of \uU\ having at least two elements.

Of course, polyreferential models are also traditional, and we use the same interpretation of formulas as in all traditional models.

We say that a formula is a \emph{polyreferential} (\emph{set-theoretic}) \emph{tautology} if and only if it is true in all polyreferential models. Since \Lcal\ is sound with respect to all traditional models, it is also sound with respect to all polyreferential models. In remarks of Section~\ref{sec6}, we will show that \Lcal\ is complete with respect to all polyreferential models.

\section{Other possible formal approaches to syllogistic\label{sec3}}

In this paper, we omit the disputes \citep[cf., e.g.,][]{Cor72a, Cor74, Smil} over what Aristotle’s syllogisms are and what form Aristotle’s syllogistic itself has, as it is irrelevant here. Below, we will show that, according to Corcoran’s views, we can treat syllogisms as argument forms and apply other types of inference rules to them. Smiley’s views on syllogisms and syllogistic will be presented at the end of this section.

We use the notations of the form $\pi_1,\ldots,\pi_n\Longrightarrow\omega$ for (valid) argument schemes, where  lowercase Greek letters represent arbitrary sentence formulas as their \pagebreak  premises and conclusion.\footnote{We assume that repeated premises and their order are unimportant in a given sequent. Therefore, we identify sequents that differ in one or both of these features.}  Moreover, the notations of the form $\Longrightarrow\omega$ are to express tautologies. A~sequent $\pi_1,\ldots,\pi_n\Longrightarrow\omega$ corresponds to the implication  $(\pi_1\wedge\dots\wedge\pi_n)\to\omega$.

\looseness=-1 One of the possible solutions is to use an appropriate natural deduction system with inference rules corresponding to acceptable argument forms and the so-called \emph{proof construction rules}. With their help, we derive new argument forms. For example, \citet{Cor72b, Cor74} proposes an understanding of Aristotelian syllogistic, in which the proof construction rule consists of assumption proofs, and selected correct syllogistic modes (as valid argument forms) are treated as rules of inference. Corcoran presented this approach in his polemic on Łukasiewicz. The presentations of Aristotle’s syllogistic in the form of a natural deduction system can be found, among others, in \citep{Moss,Clark, Sm}.

Another possible solution is to reconstruct traditional logic using sequent calculus. As axiomatic, we can accept the sequents corresponding to all substitution instances of the axioms of \Lcal\ for all $\cS, \cP, \cM\in\mathrm{GN}$ (cf.\ Remark~\ref{iLcal}):
\begin{gather*}
	\Longrightarrow\cS\sa\cS \tag*{\eqref{Ia}}\\
\Longrightarrow\cS\si\cS \tag*{\eqref{Ii}}\\
\cM\sa\cP,\cS\sa\cM\Longrightarrow\cS\sa\cP \tag*{\eqref{Barbara}}\\
\cM\sa\cP,\cM\si\cS\Longrightarrow\cS\si\cP \tag*{\eqref{Datisi}}\\
\cS\se\cP\Longleftrightarrow\neg\,\cS\si\cP \tag*{\eqref{dfe}}\\
\cS\so\cP\Longleftrightarrow\neg\,\cS\sa\cP \tag*{\eqref{dfo}}
\end{gather*}
where for formulas $\alpha$ and $\beta$, $\alpha\Longleftrightarrow\beta$ is short for two sequents $\alpha\Longrightarrow\beta$ and $\beta\Longleftrightarrow\alpha$.

Moreover, to facilitate the derivation of successive sequents, all sequents corresponding to all substitutions of all consequents and tautologies in CPL with formulas from $\ForL$ can be taken as CPL-axioms. In other words, if we have a consequence $\varphi_1,\ldots,\varphi_n\vDash\psi$ in CPL, then we take as an axiom the sequent  $\pi_1,\ldots,\pi_n\Longrightarrow\omega$, obtained from the consequent by substituting the propositional letters with formulas from $\ForL$. For example, we have axiomatic sequents obtained from $\alpha\Longleftrightarrow \neg\,\neg\,\alpha$.

Further, we will write finite sequents of premises (possibly empty) using capital Greek letters. To “generate” sequents,  we use derivation (inference) rules such as the cut rule. In the considered case, it will have the following form:
\[
\frac{A\Longrightarrow\alpha\qquad\quad B,\alpha\Longrightarrow\omega} {A,B\Longrightarrow\omega}
\]
For example, using the cut rule, from the sequents $\Longrightarrow S\si S$ and $S\sa P, S\si S\Longrightarrow S\si P$ we get the sequent $S\sa P \Longrightarrow S\si P$ corresponding to \eqref{aSi}. Applying the cut rule to the sequents $\Longrightarrow P\sa P$ and $P\sa P,P\si S\Longrightarrow S\si P$ we get the sequent $P\si S\Longrightarrow S\si P$ corresponding to \eqref{Ci}.

By the cut rule and the CPL-axioms $\alpha,\beta \Longrightarrow \alpha\wedge\beta$, $\alpha\wedge\beta \Longrightarrow \alpha$, $\alpha\wedge\beta \Longrightarrow \beta$, $\alpha,\beta \Longrightarrow \alpha\wedge\beta$, $\alpha,\alpha\to\beta \Longrightarrow \beta$, $\alpha\leftrightarrow\beta,\alpha\Longrightarrow \beta$, $\alpha\leftrightarrow\beta,\beta \Longrightarrow \alpha$ and $\alpha\to\beta , \allowbreak \beta\to\alpha\Longrightarrow\alpha\leftrightarrow\beta$ we obtain the following derivable rules:
\begin{gather*}
\frac{\varPi,\alpha,\beta\Longrightarrow\omega} {\varPi,\alpha\wedge\beta\Longrightarrow\omega}
\quad\qquad
\frac{\varPi,\alpha\wedge\beta\Longrightarrow\omega} {\varPi,\alpha,\beta\Longrightarrow\omega}
\quad\qquad
\frac{\varPi\Longrightarrow\alpha\quad \varPi'\Longrightarrow\beta} {\varPi,\varPi'\Longrightarrow\alpha\wedge\beta}\\
\frac{\varPi\Longrightarrow\alpha\to\beta}  {\varPi,\alpha\Longrightarrow\beta}\qquad\quad
\frac{\varPi\Longrightarrow\alpha\to\beta\quad\varPi\Longrightarrow\alpha}  {\varPi\Longrightarrow\beta}\\
\frac{\Longrightarrow\alpha\leftrightarrow\beta}{\alpha\Longleftrightarrow\beta}      \qquad\quad
\frac{\varPi\Longrightarrow\alpha\to\beta\quad \varPi'\Longrightarrow\beta\to\alpha} {\varPi,\varPi'\Longrightarrow\alpha\leftrightarrow\beta}
\end{gather*}
In addition, as primary, the following deduction rule must be adopted:
\[
\frac{\varPi,\alpha\Longrightarrow\omega}{\varPi\Longrightarrow\alpha\to\omega} \]
Notice that using the deduction rule, by $\alpha\to\beta,\neg\beta\Longrightarrow\neg\alpha$,  $\alpha\to\neg\beta,\beta\Longrightarrow\neg\alpha$, $\neg\alpha\to\beta,\beta\Longrightarrow\neg\alpha$ (or $\alpha\to\beta\Longleftrightarrow\neg\beta\to\neg\alpha$, $\alpha\to\neg\beta\Longleftrightarrow\beta\to\neg\alpha$ and $\neg\alpha\to\beta\Longleftrightarrow\neg\beta\to\alpha$) we obtain the following derivable contraposition rules:
\[
\frac{\varPi,\alpha\Longrightarrow\beta}{\varPi,\neg\beta\Longrightarrow\neg\alpha} \qquad \frac{\varPi,\neg\alpha\Longrightarrow\neg\beta}{\varPi,\beta\Longrightarrow\alpha} \qquad \frac{\varPi,\neg\alpha\Longrightarrow\beta}{\varPi,\neg\beta\Longrightarrow\alpha} \qquad \frac{\varPi,\alpha\Longrightarrow\neg\beta}{\varPi,\beta\Longrightarrow\neg\alpha} \]
Indeed, for example, suppose that $\varPi,\alpha\Longrightarrow\beta$. Then, by the deduction rule, we have $\varPi\Longrightarrow\alpha\to\beta$. Hence, by $\alpha\to\beta,\neg\beta\Longrightarrow\neg\alpha$ and cutting, we get $\varPi,\neg\beta\Longrightarrow\neg\alpha$. Similarly for the others rules.

Moreover, using cutting, contraposition and $\alpha\vee\beta\Longleftrightarrow\neg(\neg\alpha\wedge\neg\beta)$ we obtain the following derivable rules for disjunction:
\[
\frac{\varPi,\alpha\vee\beta\Longrightarrow\omega}{\varPi,\alpha\Longrightarrow\omega} \qquad\quad \frac{\varPi,\alpha\vee\beta\Longrightarrow\omega}{\varPi,\beta\Longrightarrow\omega} \qquad\quad \frac{\varPi,\alpha\Longrightarrow\omega\quad \varPi',\beta\Longrightarrow\omega}{\varPi,\varPi',\alpha\vee\beta\Longrightarrow\omega}
\]
Indeed, for example, suppose that $\varPi,\alpha\Longrightarrow\omega$ and $\varPi,\beta\Longrightarrow\omega$. Then, by contraposition, we get $\varPi,\neg\omega\Longrightarrow\neg\alpha$ and $\varPi,\neg\omega\Longrightarrow\neg\beta$. So, by a obtained derivable rule, we get $\varPi,\neg\omega\Longrightarrow\neg\alpha\wedge\neg\beta$ and $\varPi,\neg(\neg\alpha\wedge\neg\beta)\Longrightarrow \omega$. Hence, by $\alpha\vee\beta\Longleftrightarrow\neg(\neg\alpha\wedge\neg\beta)$ and cutting, we get $\varPi,\alpha\vee\beta\Longrightarrow\omega$. Similarly for the others rules.

Using the above rules, we obtain two new derivable rules, giving the connection between syllogisms as formulas and as sequents:
\[
\frac{\pi_1,\ldots,\pi_n\Longrightarrow\omega} {\Longrightarrow(\pi_1\wedge\dots\wedge\pi_n)\to\omega} \qquad\qquad \frac{\Longrightarrow(\pi_1\wedge\dots\wedge\pi_n)\to\omega}  {\pi_1,\ldots,\pi_n\Longrightarrow\omega}
\]
\indent The above consideration shows that the reconstruction of Aristotle’s syllogistic as a sequent calculus is equivalent to the reconstruction given by Łukasiewicz.

Another reconstruction of Aristotle’s syllogistic was given by \citet{Smil}. He writes (p.~139):
\begin{quote}\small
Given that Aristotle is concerned with deductions, i.e., with how conclusions may be derived, we should expect him to be equally concerned with deducibility, i.e., with what conclusions are derivable. We should also bear in mind that deducibility can be discussed either by means of verbs such ‘as \dots\ implies \dots’ or ‘\dots\ follows from \dots’, or by means of conditionals such as ‘if \dots\ then necessarily \dots’ or plain ‘if \dots\ then \dots’; the difference between the verbal form and the conditional form being merely the difference between mention and use. In this way think we can explain Aristotle’s frequent use of conditionals in his discussions of syllogistic without needing to identify, as Łukasiewicz does, the conditionals with the syllogisms themselves.
\end{quote}
Smiley says that Aristotle’s “argumentative structure” is suitably expressed by a \emph{proof- sequence} or \emph{deduction}. To reconstruct Smiley’s notion of \emph{deduction}, let us assume that $\pi$  and $\alpha$ represent schemas of categorical sentences, and $\varPi$ represents their sets. Standardly, for all name letters $\cS$ and $\cP$, $\cS\sa\cP$ and $\cS\so\cP$ (resp.\ $\cS\si\cP$ and $\cS\se\cP$) are mutually contradictory. \citet[p.~141]{Smil} accepts the following “rules of inference”:
\begin{enumerate}[\textbullet]
\item $\cS\sa\cM,\cM\sa\cP\rhd_1 \cS\sa\cP$
\item $\cS\sa\cM,\cM\se\cP\rhd_2 \cS\se\cP$
\item $\cP\se\cS\rhd_3 \cS\se\cP$
\item $\cP\sa\cS\rhd_4 \cS\si\cP$
\end{enumerate}
Two of them correspond to the syllogisms \emph{Barbara} and \emph{Celarent}, the third is the law of \emph{conversion} of universal denial sentences, and the fourth is the law of \emph{conversion per accidens} of affirmative sentences. Smiley reconstruction of Aristotle’s syllogistic is thus a kind of natural deduction system. We will show that it can also be represented as a kind of sequent calculus.

To signify that there exists a deduction of $\omega$ from $\varPi$, \citet[p.~142]{Smil} writes $\varPi\vdash\omega$. The definition of formal deduction Smiley gives inductively. We will write by using one trivial deduction ($\alpha\vdash\alpha$), a certain form of cutting and the rule \emph{reductio ad impossibile}:
\[
\frac{\varPi_1\vdash\alpha_1\quad\varPi_2\vdash\alpha_2\quad \alpha_1,\alpha_2\rhd_i\omega}{\varPi_1,\varPi_2\vdash\omega}         \quad\text{for $i=1,\ldots,4$}
\]
\[
\frac{\varPi_1,\overline{\omega}\vdash\alpha\quad\varPi_2\vdash \overline{\alpha}}{\varPi_1,\varPi_2\vdash\omega}
\]
where $\overline{\alpha}$ is the formula contradicting $\alpha$.\footnote{In the cut rule, Smiley used the following form for any $n$, although $n=2$ is sufficient since inference rules have one or two premises:
\[
\frac{\varPi_1\vdash\alpha_1\quad\varPi_n\vdash\alpha_n\quad \alpha_1,\ldots,\alpha_n\rhd_i\omega}{\varPi_1,\ldots,\varPi_n\vdash\omega}         \quad\text{for $i=1,\ldots,4$}
\]} Of course, using trivial deductions and cutting, for every $i=1,\ldots ,4$, we obtain:
\begin{enumerate}[\textbullet]
\item If $\alpha_1,\alpha_2 \rhd_i \omega$, then $\alpha_1,\alpha_2\vdash\omega$.
\end{enumerate}

In the above way, Smiley reconstructed Aristotle’s syllogistic (without the fourth figure). We see that his reconstruction can be represented as a sequent calculus. However, here, the sequents represent deductions, not argument forms.

Let us note, however, that Łukasiewicz rejected the alternative approaches to the reconstruction of Aristotle’s syllogistic presented above. He believed that Aristotle’s syllogisms have the form of implications.

\section{Calculi allowing empty names\label{sec4}}

We will present calculi which can also be applied to empty general names. As a standard, we will assume that both functors of affirmative sentences will be primitive. The remaining functors will be definable using the primitive ones. We understand the functor of particular affirmative sentences with its natural interpretation, and we leave the abbreviation ‘\si’ for it. The problem, however, is that the functor for universal affirmative sentences has two variations that differ when applied to empty names. As primitive, we will take the weak interpretation in the first two points of this section and the strong one in the rest.

\subsection{Shepherdson's approach\label{subsec:4.1}}
\emph{Shepherdson’s \textup{\textsf{ai}}-system.} For the weak interpretation of the functor for universal affirmative sentences, we will leave the abbreviation ‘\sa’ and---according to the proposal of Kotarbiński and Lejewski---we can read it as ‘all \dots\ is \dots’. The set of \sa\si-formulas is built in the standard way from the atomic formulas, the Boolean propositional connectives and brackets. For the set of the primitive functors ‘\sa’ and ‘\si’, \citet{Sh} proposed an axiomatisation of the \sa\si-system. Of Łukasiewicz’s four axioms, he left \eqref{Ia}, \eqref{Barbara} and \eqref{Datisi} but rejected the principle of identity \eqref{Ii} since it turns into a false sentence for all empty names. Instead, Shepherdson took two axioms weaker than \eqref{Ii}:
\begin{gather*}
S\si P\rightarrow S\si S \label{cIi}\tag{$\star$}\\
\neg\,S\si S\rightarrow S\sa P \label{nES}\tag{$\star\star$}
\end{gather*}
The first one says that every true particular affirmative sentence has a non-empty subject. The second enforces the truth of all universal affirmative sentences with empty subjects.

Notice that now, using only \eqref{Datisi}, we obtain the polysyllogisms:
\begin{gather*}
(M\si Q \wedge M\sa P \wedge Q\sa S)\to S\si P \label{Datisiplus} \tag{Datisi+}\\
(M\si M\wedge M\sa P \wedge M\sa S)\to S\si P \label{Daraptiplus} \tag{Darapti+}
\end{gather*}
Of course, from \eqref{Datisiplus} and \eqref{Ia} we obtain  \eqref{Datisi}.

\medskip

\paragraph{\emph{The Shepherdsonian \textup{\textsf{aieo}}-system.}} We can extend the set of $\sa\si$-formulas to the set of \sa\si\se\so-formulas as for Łukasiewicz's calculus. To Shepherdson's axioms, we add \eqref{dfe} and \eqref{dfo}. Moreover, all substitutions of all CPL tautologies with \sa\si\se\so-formulas are also accepted as axioms. By \Sh, we denote the system having Shephedson's axioms, definitions \eqref{dfe} and \eqref{dfo}, and two rules for deriving theses: detachment and substitution. We remember that using these means, from \eqref{Ia} and \eqref{Datisi}, we get \eqref{Ci}. Moreover, by \eqref{Ci} and \eqref{dfe}, we get \eqref{Ce}. Of course, all theses of the Shepherdsonian \sa\si\se\so-system are also theses of~\Lcal.
\begin{Remark}\label{iSh}
We can also consider another version of \Sh, which we obtain by rejecting the rule of uniform substitution and accepting all substitution instances of the axioms of \Sh\ as its specific axioms. We can show that both versions have the same theses \citep[{cf.}][]{Q50}. This second version was used, among others, in \citep{Sh}. \qed
\end{Remark}

\paragraph{\emph{Definitional extensions of} \Sh.} Let us define the unary functor ‘\ex’ (“exists”), with which we state the non-emptiness of a given name:
\[\label{dfex}\tag{df\,{\ex}}
\ex S\leftrightarrow S\si S
\]
\begin{Remark} The formula ‘$\ex S$’ is a thesis of \Lcal, with respect to its axiom \eqref{Ii}. Therefore, it makes no sense to introduce definition \eqref{dfex} in Łukasiewicz’s calculus.
\end{Remark}
From axioms \eqref{cIi} and \eqref{nES}, definition \eqref{dfex} and thesis \eqref{Ci} we have, respectively:
\begin{gather*}
S\si P\rightarrow(\ex S\wedge\ex P)\\
\neg\,\ex S\rightarrow S\sa P
\end{gather*}
From the above, \eqref{Datisi} and \eqref{dfex} we obtain:
\begin{gather*}
(S\sa P\wedge\ex S)\rightarrow S\si P\\
(S\sa P\wedge\ex S)\rightarrow\ex P
\end{gather*}
\indent Moreover, from \eqref{dfo}, \eqref{dfe} and the theses already obtained, using CPL, we get:
\begin{gather*}
(\neg\,\ex S\vee \neg\,\ex P)\rightarrow S\se P\\
S\so P\rightarrow\ex S\\
(S\se P\wedge\ex S)\rightarrow S\so P\\
(S\se P\wedge\ex P)\rightarrow P\so S\\
(\ex S\wedge\neg\,\ex P)\rightarrow S\so P
\end{gather*}
\indent For the strong interpretation of the functor of universal affirmative sentences, we adopt the abbreviation ‘\ska’ and---following Kotarbiński and Lejewski---we can read it as ‘every \dots\ is \dots’. For ‘\ska’, we adopt the following definition:
\[\label{dfka}\tag{df\,{\ska}}
S\ska P \:\leftrightarrow\: (\ex S\wedge S\sa P)
\]
So we obtain the following theses:
\begin{gather*}
S\ska P\rightarrow S\si P \label{kaSi}\tag{{\ska}S{\si}}\\
S\ska P\rightarrow(\ex S\wedge\ex P)
\end{gather*}
\indent Now let us---following Kotarbiński---introduce two symmetrical equality functors for sentences of the form ‘All $S$ is a $P$ and vice versa’ and ‘Every $S$ is a $P$ and vice versa’. Instead of those, \citet{Lej} used sentences of the form ‘$S$ is identical with $P$’ (“the functor of weak identity”) and ‘Only every $S$ is a $P$’ (“the functor of strong identity”), respectively. What we are talking about here is the identity of the extensions of two general names. Let us take the abbreviations `$\circeq$' and ‘$\doteq$’ for these functors and the definitions:
\begin{gather*}
S\mathord{\circeq} P\leftrightarrow(S\sa P\wedge P\sa S) \label{dfceq} \tag{df\,$\circeq$}\\
S\mathord{\doteq} P\leftrightarrow(S\ska P\wedge P\ska S)\label{dfdeq} \tag{df\,$\doteq$}
\end{gather*}
\indent As we remember, apart from the \emph{weak} interpretation of universal denial sentences, we have two other interpretations: \emph{strong} and \emph{super-strong}. To express the former, we can use the functor ‘no \dots\ is \dots’, leaving ‘\se’ and the definition (df e). For the strong and super-strong interpretations, we can use the functors ‘every \dots\ is not \dots’ and ‘every \dots\ is not \dots\ and vice versa’, the abbreviations ‘\ske’ and ‘\skke’, and the following definitions:%
\begin{gather*}
S\ske P\leftrightarrow(\ex S\:\wedge S\se P) \label{dfke}\tag{df\,$\ske$}\\
S\skke P\leftrightarrow(\ex S\wedge\ex P\wedge S\se P) \label{dfkke}\tag{df\,{\skke}}
\end{gather*}
So the functor `\skke' is symmetrical, but `\ske' is not. Moreover, we have:
\begin{gather*}
S\ske P\rightarrow S\so P\\
S\skke P\rightarrow(S\so P\wedge P\so S)
\end{gather*}
\paragraph{\emph{Set-theoretic semantics for} \Sh.} In the semantic study of \Sh, we can use set-theoretic semantics as for \Lcal. Now, however, we use models of the form $\la\uU,\DD\ra$, where the universe $\uU$ is an arbitrary set (may be empty) and the denotation function  $\DD$ assigns to name letters arbitrary subsets of $\uU$. With this only change, in the same way as in the first point of Subsection~\ref{subsec2.3}, we define the notions of a \emph{formula being true in a model} and of \emph{being a tautology}.

For definitional extensions of \Sh, we re-introduce the notion of \emph{formula being true in a model}. For any model $\mM=\la \uU ,\DD \ra$, for all $\cS, \cP\in\mathrm{GN}$, we assume:
\begin{enumerate}[\textbullet]
\item $\ex\cS$ is true in $\mM$ iff the set $\DD(\cS)$ is non-empty;
\item $\cS\ska\cP$ is true in $\mM$ iff the set $\DD(\cS)$ is non-empty and is included in $\DD(\cP)$;
\item $\cS\mathord{\circeq}\cP$ is true in $\mM$ iff $\DD(\cS)=\DD(\cP)$;
\item $\cS\mathord{\doteq}\cP$ is true in $\mM$ iff $\DD(\cS)=\DD(\cP)$ and the set $\DD(\cS)$ is non-empty;
\item $\cS\ske\cP$ is true in $\mM$ iff the set $\DD(\cS)$ is non-empty and is disjoint with $\DD(\cP)$;
\item $\cS\skke\cP$ is true in $\mM$ iff the sets $\DD(\cS)$ and $\DD(\cP)$ are non-empty and they are disjoint.
\end{enumerate}
For formulas built with propositional connectives, we use the standard truth tables. We say that a formula is a (\emph{set-theoretic}) \emph{tautology} (or is \emph{valid}) if and only if it is true in all models. With the above interpretation, all accepted definitions are tautologies. As can be easily seen, \Sh\ is sound with regard to set-theoretic semantics. Indeed, all axioms of \Sh\ are tautologies, uniform substitution and the detachment rule preserves the validity of formulas. In Section~\ref{sec6}, we show different ways of proving the completeness of~\Sh.
\begin{Remark}\label{rem:ne-models}
It can be assumed that model universes are non-empty. Namely,
\begin{enumerate}[\textbullet]
\item A formula is a tautology if and only if it is true in every model with a non-empty universe.
\end{enumerate}
Indeed, if $\alpha$ is true in every model having a non-empty universe, then it is true in a model with $\DD(\cS)=\emptyset$ for each $\cS\in\mathrm{GN}$. Hence, $\alpha$ is also true in the model with $\uU=\emptyset$.

Moreover, similar to footnote~\ref{footD}, the model can be just the mapping $ \DD $ itself, which assigns to any name letter an arbitrary set.\qed
\end{Remark}

\paragraph{\emph{Non-monoreferential set-theoretic semantics for Shepherdsonian system.}} We remember that in using Łukasiewicz's calculus, we can exclude monoreferential names (empty names are excluded out of necessity). In Section~\ref{sec6}, we will show that Shepherdsonian systems can apply only for non-monoreferential (i.e.\ empty or polyreferential) names. To those names correspond to \emph{non-monoreferential models} of the form $\la \uU ,\DD \ra$, where:
\begin{enumerate}[\textbullet]
\item \uU\ is a set having at least two elements,
\item \DD\ is a function of \emph{denotation}, which assigns to any name letter either the empty set or a subset of \uU\ having at least two elements.
\end{enumerate}
Of course, non-monoreferential models are also standard, and we use the same interpretation of formulas as in all models.

We say that a formula is a \emph{non-monoreferential tautology} if and only if it is true in all non-monoreferential models. Since \Sh\ is sound with respect to all models, it is also sound with respect to all non-monoreferential models. In remarks of Section~\ref{sec6}, we will show that \Sh\ is complete with respect to all non-monoreferential models.

\subsection{Słupecki's approach\label{subsec4.2}}

\emph{Słupecki's system.} \citet{Sl} proposed a calculus of names in which the functors of affirmative sentences were primary with adopted the strong interpretation for universal affirmative sentences. Therefore, we can abbreviate these functors by ‘\ska’ and ‘\si’, respectively. The theses of Słupecki’s calculus can also be applied to empty names. This system includes all correct Aristotle’s syllogisms, the laws of the logical square and conversion laws.

Słupecki adopted four \ska\si-tautologies as axioms. The first is the law of conversion \eqref{Ci}, the second is the law of  subordination \eqref{kaSi} and the others are two syllogisms:
\begin{gather*}
(M\ska P\wedge S\ska M)\rightarrow S\ska P \label{Bkarbara}\tag{B\.{a}rb\.{a}r\.{a}}\\
(M\ska P\wedge S\si M)\rightarrow S\si P \label{Dkarii}\tag{D\.{a}rii}
\end{gather*}
\indent Moreover, Słupecki also adopts \eqref{dfe} and a specific definition of the functor of particular denial sentences. It cannot be \eqref{dfo} because it has ‘\sa’, not ‘\ska’. Visually, however, the definition adopted by Słupecki corresponded to \eqref{dfo} because he used the letter ‘\sa’ but understood it in the strong sense. Since we have established the meanings of the symbols ‘\sa’, ‘\ska’ and ‘\so’, we cannot replace ‘\sa’ with ‘\ska’ in \eqref{dfo}, leaving the symbol ‘\so’. We need to replace the latter with another symbol. Let  us assume that this symbol is ‘\textsf{\~{o}}’ and that the definition adopted for it is:
\[\label{deko}\tag{df\,\textsf{\~{o}}}
S\textsf{\~{o}} P\leftrightarrow\neg\,S\ska P
\]
Moreover, all substitutions of CPL tautologies with \ska\si\se\textsf{\~{o}}-formulas are accepted as axioms. We also have two derivation rules: detachment and substitution.

Notice that accepting \eqref{deko} causes some interpretation complications. According to the adopted interpretation for ‘\ska’, we will get the interpretation of ‘\textsf{\~{o}}’, which is not consistent with the linguistic usage for particular denial sentences (cf.\ point~\ref{subsec:1.3}). Namely, it turns out that a sentence of the form ‘$S\textsf{\~{o}}P$’ is to be true iff either the name $S$ is empty or it has a referent which is not a referent of the name~$P$. Słupecki himself saw this \citep[{cf.}][p.~189]{Sl}. He, therefore, tried to circumvent the difficulty by advising that
\begin{quote}
[\dots] the sentence $Oab$ [corresponds to our `$S\textsf{\~{o}}P$'] understand only as an abbreviation of the sentence $NUab$ [corresponds to our `$\neg\,S\ska P$'] and read: it is not the case that every $a$ is a $b$.\footnote{The author of this paper translates the Polish text from \citep{Sl}.}
\end{quote}
Thus, we are to reject the original finding that his ‘$O$’ (corresponding to ‘\textsf{\~{o}}’) is the symbolic notation of the functor ‘some \dots\ is not \dots’ and assume that it is only the symbolic notation of the phrase ‘it is not the case that every \dots is \dots’. The consequences of this are as follows. The meaning of ‘$S\textsf{\~{o}}P\leftrightarrow\neg\,S\ska P$’ is just an abbreviation of the identity ‘$\neg\,S\ska P\leftrightarrow\neg\,S\ska P$’. Similarly, ‘$(P\ska M\wedge S\textsf{\~{o}}M)\to S\textsf{\~{o}}P$’ and ‘$(M\textsf{\~{o}} P \wedge M\ska S)\to S\textsf{\~{o}}P$’ are only abbreviations for ‘$(P\ska M\wedge\neg\, S\ska M)\to\neg\, S\ska P$’ and ‘$(\neg M\ska P\wedge M\ska S)\to\neg\, S\ska P$’ obtained from \eqref{Bkarbara} after substitution and the contraposition of CPL. In the alphabet of Słupecki’s calculus, no symbol would represent the functor of particular denial sentences.

As already mentioned, in Słupecki’s system, all Aristotelian syllogisms, as well as the logical square and conversion laws written with ‘\ska’, ‘\si’, ‘\se’ and ‘\textsf{\~{o}}’, obtain. However, as shown in \citep{ja87}, the theses of this system are not, for example, the following \ska\si-tautologies: \eqref{cIi} and\vspace{-6pt}
\begin{gather*}
S\si S\rightarrow S\ska S\\
S\si P\rightarrow S\ska S\label{I}\tag{\dag}\\
S\ska P\rightarrow S\si S\label{II}\tag{\ddag}\\
S\ska P\rightarrow S\ska S\\
S\ska P\rightarrow P\ska P
\end{gather*}
However, one cannot claim that Słupecki did not want to obtain implications with identities in consequents  because, by \eqref{Dkarii}, \eqref{Ci} and \eqref{kaSi}, we get:
\[\label{procent}\tag{\%}
P\ska S\rightarrow S\si S
\]
\paragraph{\emph{Complete axiomatisations of \textup{\ska\textsf{i}}-tautologies.}} In \citep{ja87, ja91b, ja91c}, it was shown that the following four sets form full axiomatisations of all \ska\si-tautologies:
\begin{enumerate}[A.]
\item Słupecki's axioms plus formula \eqref{I};
\item \eqref{Ci}, \eqref{Bkarbara}, \eqref{Dkarii} plus formulas \eqref{I} and \eqref{II};
\item \eqref{kaSi}, \eqref{Bkarbara} plus \eqref{I} and the following formula
\end{enumerate}\vspace*{-6pt}
\[\label{Dkatisi}\tag{D\.{a}tisi}
(M\ska P\wedge M\si S)\rightarrow S\si P
\]
\begin{enumerate}[A.]
\setcounter{enumi}{3}
\item \eqref{Bkarbara}, \eqref{Dkatisi}, \eqref{I} and \eqref{II}.
\end{enumerate}

To the complete axiomatisations of \ska\si-tautologies, we can add the following definition of the functor `\sa':
\[\label{dfa}\tag{df\,{\sa}}
S\sa P\leftrightarrow(\neg\,S\ska S\vee S\ska P)
\]
It gives us:\vspace{-3pt}
\[
S\sa P\leftrightarrow(\neg\,S\si S\vee S\ska P)
\]
Having `\sa', we can introduce \eqref{dfo} and definitions of other functors given in point~\ref{subsec:4.1}.

In \citep{ja91b}, the definitional equivalence of \Sh\ with the four equivalent systems for \ska\si\se\so-tautologies was demonstrated. So, these systems are complete. In \citep{ja91c}, Henkin’s method proved this.

\section{The modern syllogistic with Leśniewski’s copula\label{sec5}}

\subsection{On Leśniewski’s copula and related functors\label{subsec5.1}}

\emph{Leśniewski’s singular sentences and sentences about the identity.} The copula ‘is’ is the only primitive of Leśniewski’s Ontology. This theory can be classified as a quantifier calculus of names. In this work, however, we deal only with the quantifier-free calculus of names.

\looseness=-1 Leśniewski applied his Ontology to all names without dividing them into proper names and explicit or implicit descriptions and distinguishing whether they are general or singular. Moreover, his theory is applicable to all names: empty, monoreferential and polyreferential. He had one type of variable for all names. Since, in this paper, we are interested in Ontology only in the context of the (quantifier-free) logic of names, we will use schematic letters instead of variables. Affirmative sentences with the copula ‘is’ Leśniewski understood as follows:
\begin{enumerate}[\textbullet]
\item `$S$ is a $P$' is true if and only if the name $S$ has exactly one referent which is a referent of the name $P$.
\end{enumerate}
Leśniewski’s copula ‘is’ well be standardly symbolized by the Greek letter ‘$\sis$’ (which refers to the Latin ‘est’). So, Leśniewski's sentences have the symbolic notation ‘$S\sis P$’.

Leśniewski also used denial singular sentences of the form ‘$S$ is not a $P$’ (in short, ‘$S\bar{\sis} P$’) and sentences  about the identity of objects (in short, ‘$S\idsf P$’). His denial sentences are not equivalent to their affirmative counterparts. Leśniewski understood these two sentences as follows:
\begin{enumerate}[\textbullet]
\item ‘$S\bar{\sis} P$’ is true if and only if the name $S$ has exactly one referent which is not a referent of the name P.
\item ‘$S\idsf P$’ is true if and only if the names $S$ and $P$ have the same (one) referent.
\end{enumerate}
For `$\bar{\sis}$' and `\idsf' Leśniewski used the following definition:
\begin{align*}
S\bar{\sis} P &\leftrightarrow(S\sis S\wedge \neg\,S\sis P) \label{dfnsis}\tag{$\mathrm{df}\,{\bar{\sis}}$}\\
S\idsf P&\leftrightarrow(S\sis P\wedge P\sis S) \label{dfidsf} \tag{$\mathrm{df}\,{\idsf}$}
\end{align*}
He also used \eqref{dfceq} to definite the identity of the extensions of names.
\begin{Remark}\label{Rem5.1}
Leśniewski’s sentences ‘$S\sis P$’, ‘$S\bar{\sis} P$’ and ‘$S\idsf P$’ should be distinguished from traditional singular sentences and standard identities having one of the form: ‘$a$ is a $P$’, ‘$a$ is not a $P$’ and ‘$a$ is identical to $b$’, where  their only singular names with exactly one referent can be inserted for ‘$a$’ and ‘$b$’, and for the letter ‘$P$’, we can use any general name (cf.\ point~\ref{subsec7.1}).

Moreover, Leśniewski’s sentences ‘$S\sis P$’ and ‘$S\bar{\sis} P$’ should be distinguished from singular sentences having one of the form: ‘This $S$ is a $P$’ and ‘This $S$ is not a $P$’, where, as in \citep{CzM}, only non-empty general names can be inserted for the letters and `this~$S$' denotes a selected object from the extension of~$S$. Other solutions are also possible, but they give rise to various difficulties of interpretation (cf.\ point~\ref{subsec7.2}). \qed
\end{Remark}
\paragraph{\emph{Set-theoretic semantics for Leśniewski's functors.}} For any model $\mM=\la\uU,\DD\ra$, we extend the notion of \emph{being a true formula} for Leśniewski’s using functors. So, for all name letters $\cS$ and $\cP$, we accept:
\begin{enumerate}[\textbullet]
\item $\cS\sis\cP$ is true in \mM\ iff $\DD(\cS)$ is a singleton whose only element belongs to~$\DD(\cP)$.
\item $\cS\bar{\sis}\cP$ is true in \mM\ iff $\DD(\cS)$ is a singleton whose only element does not belong to~$\DD(\cP)$.
\item $\cS\idsf\cP$ is true in \mM\ iff $\DD(\cS)$ and $\DD(\cP)$ are identical singletons.
\end{enumerate}
{Thus, \eqref{dfnsis} and \eqref{dfidsf}  are tautologies.}

\subsection{The quantifier-free fragment of Ontology}
\text{\citet[Theorem~3.4]{Ish}} showed that the quantifier-free fragment of Ontology (in short: quantifier-free Ontology) is axiomatizable by the following three theses:
\begin{gather*}
S\sis P\rightarrow S\sis S \label{Ish1}\tag{Ish1}\\
(M\sis P\wedge S\sis M)\rightarrow S\sis P \label{Ish2}\tag{Ish2}\\
(P\sis S\wedge S\sis M)\rightarrow S\sis P \label{Ish3}\tag{Ish3}
\end{gather*}
and all substitutions of CPL tautologies with formulas of the form $\cS\sis\cP$ plus detachment and substitution rules. Of course, by \eqref{Ish1}, instead of \eqref{Ish3}, we can take:
\[\label{Ish3p}\tag{Ish3$'$}
(P\sis S\wedge S\sis S)\rightarrow S\sis P
\]
It is easy to check that \eqref{Ish1}--\eqref{Ish3} and \eqref{Ish3p}  are $\sis$-tautologies. Furthermore, Mitio \citet{T} showed that \eqref{Ish1}--\eqref{Ish3} is a complete axiomatisation of the set of \sis-tautologies.
\subsection{The fusion of Shepherdson's system with the quantifier-free  Ontology\label{subsec:5.3}}

The copula ‘\sis’ is not definable by the pair of the functors ‘\sa’ and ‘\si’. Therefore, ‘\sis’ must be added to them as a primitive functor. In \citep{ja91a, ja91c}, four complete axiomatisations of the set of \sa\si\sis-tautologies are given. They all extend the axioms of~\Sh. In each of them, we add some \sa\si\sis-tautologies.

\smallskip
\noindent I. \eqref{Ish1} and \vspace{-3pt}
\begin{gather}
S\sis P\rightarrow S\sa P\label{isSa}\\
S\sis S\rightarrow S\si S \label{isSi}\\
(S\sa M\wedge M\sis M\wedge  S\si P)\rightarrow S\sis P \label{moje}
\end{gather}
II. \eqref{Ish1}, \eqref{isSa}, \eqref{isSi} and\vspace{-3pt}
\begin{gather}
(S\sis S\wedge S\sa P)\rightarrow S\sis P\label{4}\\
(S\sis S\wedge S\si P)\rightarrow S\sa P\label{5}\\
(S\sa P\wedge S\si S\wedge P\sis P)\rightarrow S\sis S \label{6}
\end{gather}
III. \eqref{Ish1}, \eqref{isSa}, \eqref{isSi}, \eqref{6} and\vspace{-3pt}
\begin{equation}\label{7}
(S\si P\wedge S\sis S)\rightarrow S\sis P
\end{equation}
IV. \eqref{Ish1}, \eqref{isSa}, \eqref{isSi}, \eqref{7} and
\begin{equation*}
(S\sa P\wedge P\sis S)\rightarrow S\sis S
\end{equation*}
\indent Notice that we do not need to take \eqref{Ish2} and \eqref{Ish3} as axioms. Indeed, firstly, from \eqref{isSa} and \eqref{Barbara} we have $(S\sis  M \wedge M\sis P)\to S\sa P$ and $S\sis M\to S\sa M$. From \eqref{Ish1} and \eqref{isSi} we get $S\sis M\to S\si S$ and $M\sis P\to M\sis M$. From \eqref{Datisi} we get $(S\sa P\wedge S\si S)\to S\si P$. Thus, using \eqref{moje}, we get \eqref{Ish2}. Secondly, by \eqref{Ia} and \eqref{moje}, we get $(P\sis P\wedge P\si S)\to P\sis S$. By \eqref{Ish1}, \eqref{isSa} and \eqref{isSi}, we get $S\sis P\to S\sa P$ and $S\sis P\to S\si S$. Hence, by \eqref{isSa} and \eqref{Datisi}, we get $S\sis P\to S\si P$. Now, by \eqref{Ia}, \eqref{Datisi}, we get $S\si P\to P\si S$. So, we get \eqref{Ish3p}; and \eqref{Ish3}, by \eqref{Ish1}.

From \eqref{Ish1}, \eqref{isSa}, \eqref{Barbara} and \eqref{4} we obtain:
\begin{equation}\label{8}
(S\sis M\wedge M\sa P)\rightarrow S\sis P
\end{equation}
Moreover, by \eqref{Ish1}, \eqref{isSi}, \eqref{6}, \eqref{Barbara} and \eqref{4}, we get:
\begin{equation}\label{9}
(S\si S\wedge S\sa M\wedge M\sis P)\rightarrow S\sis P
\end{equation}
\indent We also add two definitions \eqref{dfe} and \eqref{dfo} plus all substitutions of all CPL tautologies with \sa\si\se\so\sis-formulas are also accepted as axioms. By $\Shis$, we denote the system having Shepherdson's axioms, the axiom of the grup~I, definitions \eqref{dfe} and \eqref{dfo}, and two rules for deriving theses: detachment and substitution.
\begin{Remark}\label{iShis}
We can also consider another version of \Shis, which we obtain by rejecting the rule of uniform substitution and accepting all substitution instances of the axioms of \Shis\ as its specific axioms. We can show that both versions have the same theses. \qed
\end{Remark}

\paragraph{\emph{Set-theoretic semantics.}} In the semantic study of \Shis, we can use set-theoretic semantics for \Sh, additionally using the interpretation for ‘\sis’. As can be easily seen, $\Shis$ is sound regarding set-theoretic semantics. Indeed, all its axioms are tautologies, and substitution and detachment preserve the validity of formulas. In Section~\ref{sec6}, we show different ways of proving the completeness of~\Shis.
\begin{Remark}
From the completeness of \Sh\ and \Shis, we obtain that \Shis\ is a conservative extension of \Sh, i.e., every formula of \Sh\ being a thesis of \Shis\ is a thesis of~\Sh. \qed
\end{Remark}
\subsection{Systems for \ska\textsf{i}\sis-tautologies\label{subsec:4.3}}

The copula ‘\sis’ is also not definable by the pair of functors ‘\ska’ and ‘\si’. In \citep{ja91c}, complete axiomatisations of \ska\si\sis-tautologies are given by adding to any of four complete axiomatisations of \ska\si-tautologies from point~\ref{subsec4.2} counterparts of formulas from the set I of \sa\si\sis-tautologies from point~\ref{subsec:5.3}, i.e., formulas \eqref{Ish1} and\vspace{-3pt}
\begin{gather}
S\sis P\rightarrow S\ska P \label{kj}\tag{$\dot{\text{\ref{isSa}}}$}\\
(S\ska M\wedge M\sis M\wedge  S\si P)\rightarrow S\sis P \label{kt}\tag{$\dot{\text{\ref{moje}}}$}
\end{gather}
Formula \eqref{isSi} is redundant by \eqref{kj} and \eqref{kaSi}. Having the definitional equivalence of Shepherdson’s system for \sa\si-tautologies with the four equivalent systems for \ska\si-tautologies, one can show the definitional equivalence of each of the given \sa\si\sis-systems with each of the given \ska\si\sis-systems \citep[see][]{ja91c}.

\medskip

\paragraph{\emph{The completeness of axiomatisations of \textup{\ska\textsf{i}\sis}-tautologies.}} The presented axiomatisations of \ska\si\sis-systems are complete since they are definitionally equivalent to given complete \sa\si\sis-systems.

\section{Methods for the completeness of calculi of names with respect to set-theoretic semantics\label{sec6}}

In this section, we present methods for obtaining the completeness of the considered calculi with respect to set-theoretic semantics.  The first one comes from \citep{Sh}. The second method consists of the appropriate direct application of Henkin's method to calculi of names. In it, we use canonical models built for maximal consistent sets in a given calculus. We give two ways of doing this.

\subsection{Proofs of the completeness of calculi by Shepherdson’s approach}
\emph{\text{For $\Sh$.}} In\linebreak \citeyearpar{Sh}, Shepherdson takes the following open first-order conditions, which correspond to axioms \eqref{Ia}, \eqref{Barbara}, \eqref{Datisi}, \eqref{cIi} and \eqref{nES}:
\begin{enumerate}[B1.]
\item\labeltext{B1}{B1} $Aaa$
\item\labeltext{B2}{B2} $(Aab \wedge Abc)\to Aac$
\item\labeltext{B3}{B3} $(Aab \wedge Iac)\to Icb$
\item\labeltext{B4}{B4} $Iab \to Iaa$
\item\labeltext{B5}{B5} $Iaa \vee Aab$
\end{enumerate}
A $\mathrm B_1$-\emph{algebra} is a relational structure $\la S, A, I\ra$, where $S$ is a non-empty set and $A$, $I$ are binary relations such that \eqref{B1}--\eqref{B5} are satisfied for all $a, b, c$ of $S$. A~$\mathrm B_1$-algebra is called a \emph{special} $\mathrm B_1$-\emph{algebra} when $S$ consists of a set of subsets of some set $V$, $Aab$ and $Iab$ are respectively the relations $a \subseteq b$, $a\cap b\neq\emptyset$ of inclusion and intersection of sets. \citet[Theorem~8]{Sh} proved that
\begin{enumerate}[(Th8)]
\item\labeltext{Th8}{Th8} \emph{Every\/ $\mathrm B_1$-algebra is epimorphic\footnote{Briefly speaking, an epimorphism (in other words, \emph{quasi-isomorphism} or \emph{onto homeomorphism}) is an isomorphism without injection. All epimorphic relational structures give the same true formulas.} to a special\/ $\mathrm B_1$-algebra.}
\end{enumerate}
\emph{The sketch of the proof}. Let $\mathfrak B=\la S,A,I\ra$ be a B$_1$-algebra. Let us call a non-empty subset $F$ of $S$ an $I$-\emph{set} when it satisfies (for all $a,b\in S$):
\begin{enumerate}[\textbullet]
\item if $a\in F$ and $Aab$, then $b\in F$;
\item if $a,b\in F$, then $Iab$.
\end{enumerate}
Let $V$ be the family of all $I$-sets, $e(a)\coloneqq \{F\in V : a\in F\}$ and $S'\coloneqq\{e(a): a\in S\}$.

In the proof, \citet{Sh} appealed to an analogous theorem for certain richer A$_1$-algebras, for which $V$ had to consist of maximal $I$-sets. So, Shepherdson used Zorn's lemma. In \citep{ja92}, it was shown that this is not necessary for B$_1$-algebras. Namely, it suffices to note that, by \eqref{Barbara}, \eqref{Datisiplus}, \eqref{Daraptiplus}, \eqref{cIi} and \eqref{Ci}, for all $a, b\in S$, we get:
\begin{enumerate}[(\$\$)]
\item[(\$)]\labeltext{\$}{a,b} If $Iab$ holds, then $[a,b]\coloneqq\{c\in S : Aac$ or $Abc\}$ is an $I$-set of $e(a)$ and $e(b)$.
\item[(\$\$)]\labeltext{\$\$}{a,a}  $Iaa$ holds iff $[a,a]$ is an $I$-set of $e(a)$ iff $e(a)\neq\emptyset$.
\end{enumerate}
Using the definitions of B$_1$-algebras, $I$-sets and the function $e$ we obtain that:
\begin{enumerate}[\textbullet]
\item $Aab$ holds iff $e(a)\subseteq e(b)$,
\item $Iab$ holds iff $e(a)\cap e(b)\neq\emptyset$.
\end{enumerate}
So, $e\colon S\to S'$ is an epimorphism.  \qed

By the \emph{free-variable calculus} B$_1$, Shepherdson understands the formal theory ``obtained by incorporating axioms'' \eqref{B1}--\eqref{B5} ``into the propositional calculus.'' This system contains (free) individual variables, atomic formulas $Auv$ and $Iuv$ (where $u$, $v$ are variables), and formulas built up from atomic formulas by means of the propositional connectives. The axioms are all substitution instances of  \eqref{B1}--\eqref{B5} and of the axioms of propositional calculus. The rule of inference is detachment. Shepherdson writes that by Henkin's method is obtained the completeness for the free-variable calculus~B$_1$:
\begin{enumerate}[\textbullet]
\item \emph{A formula is a theorem of\/ $\mathrm{B}_1$ if and only if it is true in all\/ $\mathrm{B}_1$-algebras.}
\end{enumerate}
From the above and \eqref{Th8}, \citet[Theorem~10]{Sh} proved that
\begin{enumerate}[(Th10)]
\item\labeltext{Th10}{Th10} \emph{A formula is a theorem of\/ $\mathrm{B}_1$ if and only if it is true in all special\/ $\mathrm{B}_1$-algebras.}
\end{enumerate}

We can definitionally extend $\mathrm{B}_1$ to the free-variable calculus $\mathrm{B}_1^{\mathrm{d}}$ by adding atomic formulas $Euv$ and $Ouv$ (where $u$, $v$ are variables) with axioms being all substitution instances of \eqref{B1}--\eqref{B5} and
\begin{enumerate}[D1.]
\item\labeltext{D1}{D1} $Eab \leftrightarrow \neg Iab$
\item\labeltext{D2}{D2} $Oab \leftrightarrow \neg Aab$
\end{enumerate}
For this extension we create special $\mathrm{B}_1^{\mathrm{d}}$-algebras assuming that $Eab$ and $Oab$ are respectively the relations $a\cap b = \emptyset$, $a\nsubseteq b$. Since \eqref{D1} and \eqref{D2} are true in every special $\mathrm{B}_1^{\mathrm{d}}$-algebra, we obtain a counterpart of \eqref{Th8}. Moreover, as for $\mathrm{B}_1$, we obtain the completeness for B$_1^{\mathrm{d}}$. Hence, we obtain the counterpart of \eqref{Th10}:
\begin{enumerate}[(Th10$'$)]
\item\labeltext{Th10$'$}{Th10'} \emph{A formula is a theorem of\/ $\mathrm{B}_1^{\mathrm{d}}$ if and only if it is true in all special\/ $\mathrm{B}_1^{\mathrm{d}}$-algebras.}
\end{enumerate}

The free-variable calculus B$_1^{\mathrm{d}}$ can, of course, be identified with the version of the calculus \Sh\ in which we accept all substitution instances of the axioms of \Sh\ as its specific axioms and use only the rule of detachment (cf.\ Remark~\ref{iSh}). For every formula $\alpha$ of \Sh, let $\alpha^{\star}$ be its counterpart in the language of B$_1^{\mathrm{d}}$. So, we get
\begin{enumerate}[\textbullet]
\item $\alpha$ is a thesis of \Sh\ iff $\alpha^{\star}$ is a theorem of $\mathrm{B}_1^{\mathrm{d}}$.
\item $\alpha$ is a tautology iff $\alpha^{\star}$ is true in all special\/ $\mathrm{B}_1^{\mathrm{d}}$-algebras.
\item $\alpha$ is a thesis of \Sh\ iff $\alpha$ is a tautology.
\end{enumerate}

\smallskip

\paragraph{\emph{For} \Lcal.} In \citeyearpar{Sh}, Shepherdson creates the free-variable calculus B$_3$ by replacing \eqref{B4} and \eqref{B5} with one axiom which is a counterpart of \eqref{Ii}:
\begin{enumerate}[B4$'$.]
\item\labeltext{B4$'$}{B4'} $Iaa$
\end{enumerate}
For B$_3$, Shepherdson created special B$_3$-algebras, which differ from special B$_1$-algebras only in that their universes consist of non-empty sets since it has to meet~\eqref{B4'}. He obtains the following theorem~11:
\begin{enumerate}[(Th11)]
\item\labeltext{Th11}{Th11} \emph{Every\/ $\mathrm{B}_3$-algebra is epimorphic to a special\/ $\mathrm{B}_3$-algebra. A formula of\/ $\mathrm{B}_3$ is a theorem if and only if it is true in all special\/ $\mathrm{B}_3$.}
\end{enumerate}
Indeed, in the second variant of the proof of \eqref{Th8}, one can notice that by \eqref{B4'} and \eqref{a,a}, for each $a\in S$ we have $e(a)\neq\emptyset$.
\begin{Remark}
To Theorem~11, Shepherdson adds footnote~9, which says that a similar theorem can be obtained for Słupecki's system from \citep{Sl}, which we discussed in Section~\ref{subsec4.2}. However, as we showed there, this system is not complete. (Shepherdson stated that Słupecki's works   ``were not available to the author''; he knew them only from their reviews.)\qed
\end{Remark}

The rest is as for \Sh, only changing the tautologies to traditional tautologies and the special B$_1$-algebras to the special B$_3$-algebras. So, we obtain:
\begin{enumerate}[\textbullet]
\item $\alpha$ is a thesis of \Lcal\ iff $\alpha$ is a traditional tautology.
\end{enumerate}

\smallskip
\paragraph{\emph{For} \Shis.} For conditions \eqref{B1}--\eqref{B5} we add the following, which correspond to axioms \eqref{Ish1}, \eqref{isSa}--\eqref{moje}:
\begin{enumerate}[C1.]
\setcounter{enumi}{-1}
\item\labeltext{C0}{C0} $\varepsilon ab\to \varepsilon aa$
\item\labeltext{C1}{C1} $\varepsilon ab\to Aab$
\item\labeltext{C2}{C2} $\varepsilon aa\to Iaa$
\item\labeltext{C4}{C4} $(Aac\wedge \varepsilon  cc\wedge  Iab)\rightarrow \varepsilon ab$
\end{enumerate}
Analogously to \citep{Sh}, a C-\emph{algebra} is a relational structure $\la S, A, I,\varepsilon\ra$, where $S$ is a non-empty set and $A$, $I$ and $\varepsilon$ are binary relations such that \eqref{B1}--\eqref{B5}, \eqref{C0}--\eqref{C4} are satisfied for all $a, b, c$ of $S$. A~C-algebra is called a \emph{special} C-\emph{algebra} when $S$ consists of a set of subsets of some set $V$, $Aab$, $Iab$, $\varepsilon ab$ are respectively the relations $a\subseteq b$, $a\cap b\neq\emptyset$, $a$ is a singleton and $a\subseteq b$. We get \citep[{cf.}][]{ja92}:
\begin{enumerate}[(Th1$\varepsilon$)]
\item\labeltext{Th1$\varepsilon$}{Th1epsilon} \emph{Every\/ $\mathrm C$-algebra is epimorphic to a special\/ $\mathrm C$-algebra.}
\end{enumerate}
\emph{The sketch of the proof}. Let $\mathfrak C=\la S,A,I,\varepsilon\ra$ be a C-algebra. As for B$_1$-algebras we define $I$-sets and $I$-sets of the form $[a,b]$. Let $V$ be the sum of the family of all $I$-sets and the set~$S$. Moreover, for every $a\in S$ we put:
\[
e(a)\coloneqq
\begin{cases}
\{F\in V : a\in F\} & \quad\text{if $\varepsilon aa$,}\\
\{F\in V : a\in F\}\cup \{c\in S : Icc\text{~and~}Aca\} & \quad\text{otherwise.}\\
\end{cases}
\]
We obtain:
\begin{enumerate}[\textbullet]
\item if non $Iaa$, then $e(a)=\emptyset$;
\item if non $\varepsilon aa$ and $Iaa$, then $\{a,[a,a]\}\subseteq e(a)$;
\item if $\varepsilon aa$, then $e(a)=\{[a,a]\}$.
\end{enumerate}
Also as for B$_1$-algebras we put $S'\coloneqq\{e(a): a\in S\}$. Using the definitions of C-algebras, $I$-sets and the function $e$ we obtain:
\begin{enumerate}[\textbullet]
\item $Aab$ holds iff $e(a)\subseteq e(b)$,
\item $Iab$ holds iff $e(a)\cap e(b)\neq\emptyset$.
\item $\varepsilon ab$ holds iff $e(a)$ is a singleton and $e(a)\subseteq e(b)$.
\end{enumerate}
So, $e\colon S\to S'$ is an epimorphism.  \qed

However, instead of considering free-variable calculus, we can consider the standard open first-order theory C obtained by incorporating axioms \eqref{B1}--\eqref{B5}, \eqref{C0}--\eqref{C4}. This system contains atomic formulas $Auv$, $Iuv$ and $\varepsilon uv$ (where $u$, $v$ are individual variables). The specific axioms are substitution instances of \eqref{B1}--\eqref{B5}, \eqref{C0}--\eqref{C4} by $a/x$, $b/y$ and $c/x$. Moreover, we use the standard axioms for first-order theories as logical axioms.

As for all first-order theories, also for C, we can use Gödel’s completeness theorem. So, by \eqref{Th1epsilon}, we have:
\begin{enumerate}[(Th2$\varepsilon$)]
\item\labeltext{Th2$\varepsilon$}{Th2epsilon} \emph{A formula is a theorem of\/ $\mathrm{C}$ if and only if it is true in all special\/ $\mathrm{C}$-algebras.}
\end{enumerate}
Moreover, since C is open, we consider the quantifier-free theory  $\mathrm{C}^{\mathrm{o}}$ such that:
\begin{enumerate}[($\mathrm{C}^{\mathrm{o}}$)]
\item a formula is a theorem of $\mathrm{C}^{\mathrm{o}}$ if and only if it is open and derivable from all open axioms of \textup{C} by detachment and substitution.
\end{enumerate}
By the known fact \citep[see, e.g.,][p.~329]{RiS}, we have:
\begin{enumerate}[(Th3$\varepsilon$)]
\item\labeltext{Th3$\varepsilon$}{Th3epsilon} \emph{For any open formula\/\textup{:} it is a theorem of\/ $\mathrm{C}^{\mathrm{o}}$ if and only if it is theorem of\/ \textup{C}.}
\end{enumerate}
From \eqref{Th2epsilon} and \eqref{Th3epsilon} we get:
\begin{enumerate}[(Th4$\varepsilon$)]
\item\labeltext{Th4$\varepsilon$}{Th4epsilon} \emph{For any open formula\/\textup{:} it is a theorem of\/ $\mathrm{C}^{\mathrm{o}}$ if and only if it is true in all special\/ $\mathrm{C}$-algebras.}
\end{enumerate}
Next, we continue as for B$_1$, creating the quantifier-free definitional extension C$^{\mathrm{od}}$ of the theory C$^{\mathrm{o}}$. We obtain:
\begin{enumerate}[(Th4$'\varepsilon$)]
\item\labeltext{Th4$'\varepsilon$}{Th4pepsilon} \emph{For any open formula\/\textup{:} it is a theorem of\/ $\mathrm{C}^{\mathrm{od}}$ if and only if it is true in all special\/ $\mathrm{C}^{\mathrm{d}}$-algebras.}
\end{enumerate}

The quantifier-free theory C$^{\mathrm{od}}$ can be identified with the calculus \Shis. For every formula $\alpha$ of \Shis, let $\alpha^{\star}$ be its counterpart in the language of C$^{\mathrm{od}}$. So, we get:
\begin{enumerate}[\textbullet]
\item $\alpha$ is a thesis of \Shis\ iff $\alpha^{\star}$ is a theorem of $\mathrm{C}^{\mathrm{od}}$.
\item $\alpha$ is a tautology iff $\alpha^{\star}$ is true in all special\/ $\mathrm{C}^{\mathrm{d}}$-algebras.
\item $\alpha$ is a thesis of \Shis\ iff $\alpha$ is a tautology.
\end{enumerate}
\begin{Remark}
Defining $V$ for $\mathrm{B}_1$-algebras and $\mathrm{B}_3$ as for C-algebras and defining an endomorphism $e$ for these first algebras using the second of the alternative conditions used for $\mathrm{C}$-algebras gives a family $S'$ without singletons. Hence, we get that for $\Sh$ and $\Lcal$, we can use non-monoreferential and polyreferential semantics, respectively.

Non-monoreferential semantics cannot be applied to $\Shis$ since all formulas of the form $\cS\sis\cP\to\alpha$ are non-monoreferential tautologies but not all of them are theses of~$\Shis$ (e.g., $S\sis S \to M\sis M$).  \qed
\end{Remark}

\subsection{Proofs of the completeness of calculi by direct use of Henkin's method}
In \citep{ja91a, ja91c}, it is proved the completeness of the considered calculi by direct use of Henkin's method, in which we use canonical models built for maximal consistent sets in a given calculus. For every maximal consistent set $\varGamma$ of formulas in a given calculus, we will construct an appropriate canonical model $\mM_{\varGamma} =\la \uU_{\varGamma}, \DD_{\varGamma}\ra$ in which all formulas from $\varGamma$ are true. We have two ways of doing this, with and without filters.

\subsubsection{With using filters designated by maximal consistent sets in a given calculus}

In \citep{ja91a}, a universe $\uU_{\varGamma}$ of $\mM_{\varGamma}$ consists of filters built from name letters. These filters are counterparts of the $I$-sets used by \citet{Sh}. Namely, for each of the calculi $\Sh$, $\Lcal$ and $\Shis$, for an arbitrary its maximal consistent set $\varGamma$ of formulas, a \emph{filter designated by} $\varGamma$ is a non-empty subset $\nabla$ of GN satisfying the following conditions for all $\cS,\cP\in\mathrm{GM}$:
\begin{enumerate}[\textbullet]
\item if $\cS\in\nabla$ and $\cS\sa\cP\in\varGamma$, then $\cP\in\nabla$;
\item if $\cS,\cP\in\nabla$, then $\cS\si\cP\in\varGamma$.
\end{enumerate}
Similar to the proof of \eqref{Th8}, using \eqref{Barbara}, \eqref{Datisiplus}, \eqref{Daraptiplus}, \eqref{cIi} and \eqref{Ci}, we obtain the counterparts of conditions \eqref{a,b} and \eqref{a,a} for all $\cS,\cP\in\mathrm{GN}$:
\begin{enumerate}[(f1)\enspace]
\item[(f1)\enspace]\labeltext{f1}{f1} if  $\cS\si\cP\in\varGamma$, then $[\cS,\cP]\coloneqq\{\cM : \cS\sa\cM\in\varGamma$ or $\cP\sa\cM\in\varGamma\}$ is a filter;
\item[(f2)\enspace]\labeltext{f2}{f2} if $\cS\si\cS\in\varGamma$, then  $[\cS]\coloneqq[\cS,\cS]\coloneqq\{\cM : \cS\sa\cM\in\varGamma\}$ is a filter.
\end{enumerate}
After these changes, the set $\uU_{\varGamma}$ and the function $\DD_{\varGamma}$ are defined analogously to the set $V$ and the epimorphism $e$ in the proofs of \eqref{Th8} and \eqref{Th1epsilon}, respecyively.

\medskip

\paragraph{\emph{For} \Sh.} For any maximal consistent set $\varGamma$ in $\Sh$, we use  $\mM_{\varGamma}=\langle\uU_{\varGamma}, \DD_{\varGamma}\rangle$, where:
\begin{enumerate}[\textbullet]
\item $\uU_{\varGamma}$ consists of all filters designated by $\varGamma$ ($\uU_{\varGamma}$ may be empty),
\item $\DD_{\varGamma}(\cS)\coloneqq\{\nabla\in\uU_{\varGamma} : \cS\in\nabla\}$.
\end{enumerate}
For every name letter $\cS$, we obtain:
\begin{enumerate}[(a)]
\item\label{a} $\cS\si\cS\in\varGamma$ iff   $\cS\in[\cS]\in\DD_{\varGamma}(\cS)$ iff  $\DD_{\varGamma}(\cS)\neq\emptyset$.
\end{enumerate}
Indeed, by  \eqref{f2} and \eqref{Ia}, if $\cS\si\cS\in\varGamma$, then $[\cS]$ is a filter and $\cS\in[\cS]\in\DD_{\varGamma}(\cS)$. So, $\DD_{\varGamma}(\cS)\neq\emptyset$. If  $\DD_{\varGamma}(\cS)\neq\emptyset$, then for some $\nabla$ we have $\cS\in\nabla$. So, $\cS\si\cS\in\varGamma$.

By induction, for any formula $\alpha$, we have:
\begin{enumerate}[(C\Sh)]
\item\labeltext{C\Sh}{CSh} $\alpha$ is true in $\mM_{\varGamma}$ iff $\alpha\in\varGamma$.
\end{enumerate}
\emph{The sketch of the proof}. Firstly, for ${\circ}\in\{\sa,\si,\se,\so\}$ and all $\cS,\cP\in\mathrm{GM}$, we have:
\begin{enumerate}[\textbullet]
\item $\cS\circ \cP$ is true in $\mM_{\varGamma}$ iff $\cS\circ\cP\in\varGamma$.
\end{enumerate}

Let $\cS\sa \cP$ is true in $\mM_{\varGamma}$, i.e.\ $\DD_{\varGamma}(\cS)\subseteq\DD_{\varGamma}(\cP)$. Then, if $\DD_{\varGamma}(\cS)=\emptyset$, then $\cS\si\cS\notin\varGamma$, by~\eqref{a}. Hence $\cS\sa\cP\in\varGamma$, by \eqref{nES}. If $\DD_{\varGamma}(\cS)\neq\emptyset$, then $\cS\si\cS\in\varGamma$ and  $[\cS]\in\DD_{\varGamma}(\cS)$, by \eqref{a}. So, also $[\cS]\in\DD_{\varGamma}(\cP)$. Hence $\cS\sa\cP\in\varGamma$. We have the converse implication from definitions of filters and $\DD_{\varGamma}$.

Let $\cS\si\cP$ be true in $\mM_{\varGamma}$, i.e.\ some $\nabla$ belongs to $\DD_{\varGamma}(\cS)$ and $\DD_{\varGamma}(\cP)$. Then $\cS, \cP\in\nabla$. Hence $\cS\si\cP\in\varGamma$. Conversely, assume that $\cS\si\cP\in\varGamma$. Then, by \eqref{f1} and \eqref{Ia}, $[\cS,\cP]$ is a filter for which $\cS$ and $\cP$ belong. So, $[\cS,\cP]$  belongs to $\DD_{\varGamma}(\cS)$ and $\DD_{\varGamma}(\cP)$.

For `\so' and `\se': By the above facts, using \eqref{dfo} and \eqref{dfe}, respectively.

Secondly, we can use standard properties of the connectives in maximal consistent sets and obtain that a formula is true in $\mM_{\varGamma}$ iff it belongs to~$\varGamma$.   \qed

Therefore, by \eqref{CSh}, using the properties of maximal consistent sets, we obtain:
\begin{enumerate}[\textbullet]
\item \emph{A formula is a thesis of\/ \textup{\Sh} if and only if it is an\/ \sa\si\se\so-tautology}.
\end{enumerate}
\emph{The sketch of the proof}. Suppose that $\alpha$ is a \sa\si\se\so-tautology and $\varGamma$ is an arbitrary maximal consistent set of formulas. Then, by \eqref{CSh}, $\alpha\in\varGamma$ since $\alpha$ is true in~$\mM_{\varGamma}$. So, $\alpha$ belongs to all maximal consistent sets for $\Sh$. Hence, $\alpha$ is a thesis of~$\Sh$. \qed

\medskip

\paragraph{\emph{For} \Lcal.} For any maximal consistent set $\varGamma$ in $\Lcal$, we use $\mM_{\varGamma}=\langle\uU_{\varGamma}, \DD_{\varGamma}\rangle$ as for $\Sh$. Now, by \eqref{Ii} and \eqref{a}, $\mM_{\varGamma}$ is traditional models since for each $\cS\in\mathrm{GN}$, we have that $[\cS]$ is a filter belonging to $\DD_{\varGamma}$. The rest of the proof is similar to that for \Sh. Therefore, we get:
\begin{enumerate}[\textbullet]
\item \emph{A formula is a thesis of\/ \textup{\Lcal} if and only if it is a traditional\/ \sa\si\se\so-tautology}.
\end{enumerate}

\medskip

\paragraph{\emph{For} $\Shis$.} For any maximal consistent set $\varGamma$ in $\Shis$, we use  $\mM_{\varGamma}=\langle\uU_{\varGamma}, \DD_{\varGamma}\rangle$, where:
\begin{enumerate}[\textbullet]
\item $\uU_{\varGamma}$ consists of all filters designated by $\varGamma$ and all name letters $\cS$ such that $\cS\si\cS\in\varGamma$,
\item $\DD_{\varGamma}(\cS)\coloneqq
\begin{cases}
\{\nabla\in\uU_{\varGamma} : \cS\in\nabla\} &\text{if $\cS\sis\cS\in\varGamma$}\\
\{\nabla\in \uU_{\varGamma} : \cS\in\nabla\}\cup\{\cM\in\uU_{\varGamma} : \cM\sa\cS\in\varGamma\}&\text{if $\cS\sis\cS\notin\varGamma$}
\end{cases}$
\end{enumerate}
For all name letter $\cS$ and $\cM$, we obtain:
\begin{enumerate}[(a)]
\item\label{ab} If $\cM\in\DD_{\varGamma}(\cS)$, then $\cS\sis\cS\notin\varGamma$ and $\cS\si\cS\in\varGamma$.
\item\label{bb} $\cS\si\cS\in\varGamma$ iff   $\cS\in[\cS]\in\DD_{\varGamma}(\cS)$ iff  $\DD_{\varGamma}(\cS)\neq\emptyset$.
\item\label{cb} if $\cS\sis\cS\notin\varGamma$ and $\cS\si\cS\in\varGamma$, then $\{\cS,[\cS]\}\subseteq\DD_{\varGamma}(\cS)$;
\item\label{dc} $\cS\sis\cS\in\varGamma$ iff  $\DD_{\varGamma}(\cS)=\{[\cS]\}$ iff $\DD_{\varGamma}(\cS)$ is a singleton.
\end{enumerate}
\emph{Proof}. For \eqref{ab}: By \eqref{Daraptiplus} and the definition of $\DD_{\varGamma}$. For \eqref{bb}: As for \Sh\ using \eqref{ab}. For \eqref{cb}: By \eqref{Ia} and \eqref{bb}.

For \eqref{dc}: Let $\cS\sis\cS\in\varGamma$. Then, by \eqref{isSa},  \eqref{isSi}  and \eqref{bb}, $\cS\sa\cS\in\varGamma$, $\cS\si\cS\in\varGamma$, $[\cS]$ is a filter and $\cS\in[\cS]$. We show that if $\nabla\in\DD_{\varGamma}(\cS)$ then $\nabla=[\cS]$. Let $\cM\in\nabla$. Then $\cS\si\cM\in\varGamma$. Hence, by \eqref{7}, $\cS\sa\cM\in\varGamma$. So $\cM\in[\cS]$. Conversely, if $\cM\in[\cS]$, then $\cS\sa\cM\in\varGamma$. So, $\cM\in\nabla$ since $\cS\in\nabla$. Finally, if $\DD_{\varGamma}(\cS)$ is a singleton, then $\cS\sis\cS\in\varGamma$, by \eqref{bb} and \eqref{cb}. \qed

By induction, for any formula $\alpha$, we obtain:
\begin{enumerate}[(C\Shis)]
\item\labeltext{C\Shis}{CShis} $\alpha$ is true in $\mM_{\varGamma}$ iff $\alpha\in\varGamma$.
\end{enumerate}
\emph{The sketch of the proof}. For ${\circ}\in\{\sa,\si,\se,\so,\sis\}$ and all $\cS,\cP\in\mathrm{GM}$, we have:
\begin{enumerate}[\textbullet]
\item $\cS\circ\cP$ is true in $\mM_{\varGamma}$ iff $\cS\circ\cP\in\varGamma$.
\end{enumerate}

As for \Sh, we obtain: if $\cS\sa \cP$ is true in $\mM_{\varGamma}$, then   $\cS\sa\cP\in\varGamma$. Conversely, assume that $\cS\sa\cP\in\varGamma$. If $\nabla\in\DD_{\varGamma}(\cS)$, then $\nabla\in\DD_{\varGamma}(\cP)$, by the definitions of filters and $\DD_{\varGamma}$. If $\cM\in\DD_{\varGamma}(\cS)$, then $\cM\in\DD_{\varGamma}(\cP)$, by  \eqref{Barbara}.

Let $\cS\si\cP$ be true in $\mM_{\varGamma}$, i.e.\ some $\cM$ or $\nabla$ belongs to $\DD_{\varGamma}(\cS)$ and $\DD_{\varGamma}(\cP)$. In the first case, $(\cM\si\cM\wedge\cM\sa\cS\wedge\cM\sa\cP)\in\varGamma$. So $\cS\si\cP\in\varGamma$, by \eqref{Daraptiplus}. In the second case, the same as for \Sh. The proof of the converse implication is the same as for \Sh.

Let $\cS\sis\cP$ be true in $\mM_{\varGamma}$, i.e., $\DD_{\varGamma}(\cS)$ be a singleton whose only element belongs to $\DD_{\varGamma}(\cP)$. Then, $\DD_{\varGamma}(\cS\subseteq\DD_{\varGamma}(\cP)$, and so   $\cS\sa\cP\in\varGamma$. Moreover, $\cS\sis\cS\in\varGamma$, by \eqref{dc}. Hence $\cS\sis\cP\in\varGamma$, by \eqref{4}. For the proof of the converse implication, let $\cS\sis\cP\in\varGamma$. Then $\cS\sis\cS\in\varGamma$, by \eqref{Ish1}; and so $\DD_{\varGamma}(\cS)=\{[\cS]\}$, by \eqref{dc}. Moreover, $\cS\sa\cP\in\varGamma$ and $\cS\si\cS\in\varGamma$, by \eqref{isSa} and \eqref{isSi}. Hence $\cP\in[\cS]\in\DD_{\varGamma}(\cP)$. So, $\cS\sis\cP$ is true in~$\mM_{\varGamma}$.

The rest of the proof is similar to that for \Sh. \qed

Therefore, by \eqref{CShis}, using the properties of maximal consistent sets, we get:
\begin{enumerate}[\textbullet]
\item \emph{A formula is a thesis of\/ \textup{\Shis} if and only if it is an\/ \sa\si\se\so\sis-tautology}.
\end{enumerate}
\begin{Remark}
Defining $\uU_{\varGamma}$ for $\Sh$ and $\Lcal$ as for $\Shis$ and defining $\DD_{\varGamma}$ for these first calculi using the second of the alternative conditions used for $\Shis$ gives $\DD_{\varGamma}$ with the set of values without singletons. Hence, also directly using Henkin's method, we get that for \Sh, we can use non-monoreferential semantics, and for \Lcal\ polyreferential.\qed
\end{Remark}
\subsubsection{Without using filters}

In \citep{ja91c}, the universes of models consist of simpler elements than filters designated by maximal consistent sets in a given calculus. In the case of \Sh\ and \Lcal, pairs and singletons of name letters are sufficient. In the case of \Shis, to pairs and singletons of name letters, we add the name letters themselves and their equivalence classes.

\medskip

\paragraph{\emph{For} \Sh.} For any maximal consistent set $\varGamma$ in $\Sh$, we use  $\mM_{\varGamma}=\langle\uU_{\varGamma}, \DD_{\varGamma}\rangle$, where:
\begin{enumerate}[\textbullet]
\item $\uU_{\varGamma}$ consists of all pairs $\{\cM,\cQ\}$ of name letters such that $\cM\si\cQ\in\varGamma$ ($\uU_{\varGamma}$ may be empty),
\item $\DD_{\varGamma}(\cS)$ consists of all those and only those  $\{\cM,\cQ\}{\in}\uU_{\varGamma}$ for which \mbox{$\cM\sa\cS\!\vee\!\cQ\sa\cS\in\varGamma$.}
\end{enumerate}
Of course, all singletons $\{\cM\}$ such that $\cM\si\cM\in\varGamma$ belong to $\uU_{\varGamma}$, and if also $\cM\sa\cS\in\varGamma$, then they belong to $\DD(\cS)$ (if there are any).

By  \eqref{Ia}, \eqref{cIi}, \eqref{Ci} and \eqref{Daraptiplus}, for every name letter $\cS$, we obtain:
\begin{enumerate}[(o)]
\item\labeltext{o}{o} $\cS\si\cS\in\varGamma$ iff   $\{\cS\}\in\DD_{\varGamma}(\cS)$ iff   $\DD_{\varGamma}(\cS)\neq\emptyset$.
\end{enumerate}

By induction, for any formula $\alpha$, we have:
\begin{enumerate}[(C\Sh)]
\item $\alpha$ is true in $\mM_{\varGamma}$ iff $\alpha\in\varGamma$.
\end{enumerate}
\emph{The sketch of the proof}. For ${\circ}\in\{\sa,\si,\se,\so\}$ and all $\cS,\cP\in\mathrm{GM}$, we have:
\begin{enumerate}[\textbullet]
\item $\cS\circ \cP$ is true in $\mM_{\varGamma}$ iff $\cS\circ\cP\in\varGamma$.
\end{enumerate}

The same as for filters, only changing $[\cS]$ to $\{\cS\}$ we have: if $\cS\sa\cP$ be true in $\mM_{\varGamma}$ then $\cS\sa\cP\in\varGamma$. For the proof of the converse implication, let $\cS\sa\cP\in\varGamma$ and $\{\cM,\cQ\}\in\DD(\cS)$. Then, by \eqref{Barbara}, also $\{\cM,\cQ\}\in\DD(\cP)$.

Let $\cS\si\cP$ be true in $\mM_{\varGamma}$, i.e.\ some $\{\cM,\cQ\}$ belongs to $\DD_{\varGamma}(\cS)$ and $\DD_{\varGamma}(\cP)$. Then $\cM\si\cQ\in\varGamma$, $(\cM\sa\cS\vee\cQ\sa\cS)\in\varGamma$ and $(\cM\sa\cP\vee\cQ\sa\cP)\in\varGamma$. Hence $\cS\si\cP\in\varGamma$, by \eqref{cIi}, \eqref{Ci}, \eqref{Datisiplus} and  \eqref{Daraptiplus}. For the proof of the converse implication, let $\cS\si\cP\in\varGamma$. Then, by \eqref{Ia}, $\{\cS,\cP\}$ belongs to $\DD_{\varGamma}(\cS)$ and $\DD_{\varGamma}(\cP)$.

The rest of the proof is similar to that of using filters. \qed

Therefore, by \eqref{CSh}, using the properties of maximal consistent sets,  we get:
\begin{enumerate}[\textbullet]
\item \emph{A formula is a thesis of\/ \textup{\Sh} if and only if it is an\/ \sa\si\se\so-tautology}.
\end{enumerate}
\begin{Remark}
In Remark~\ref{rem:ne-models} we show that a formula is a tautology if and only if it is true in every model with a non-empty universe. So:
\begin{enumerate}[\textbullet]
\item \emph{A formula is a thesis of\/ \Sh\ if and only if it is true in every model with a non-empty universe.}
\end{enumerate}
It can also be shown by taking canonical models with non-empty $\uU_{\varGamma}$:
\begin{enumerate}[\textbullet]
\item $\uU_{\varGamma}$ consists of all pairs of name letters,
\item $\DD_{\varGamma}(\cS)$ consists of all those and only those  $\{\cM,\cQ\}$ for which both $\cM\si\cQ\in\varGamma$ and \mbox{$\cM\sa\cS\vee\cQ\sa\cS\in\varGamma$.}
\end{enumerate}
As for the previous model, we show that condition \eqref{CSh} holds.\qed
\end{Remark}

\paragraph{\emph{For} \Lcal.} For any maximal consistent set $\varGamma$ in $\Lcal$, we use $\mM_{\varGamma}=\langle\uU_{\varGamma}, \DD_{\varGamma}\rangle$ as for $\Sh$. Now, by \eqref{Ii} and \eqref{o}, $\mM_{\varGamma}$ is traditional models. The rest of the proof is similar to that for \Sh. Therefore, we get:
\begin{enumerate}[\textbullet]
\item \emph{A formula is a thesis of\/ \textup{\Lcal} if and only if it is a traditional\/ \sa\si\se\so-tautology}.
\end{enumerate}

\smallskip

\paragraph{\emph{For} $\Shis$.} Let $\varGamma$ be a maximal consistent set in $\Shis$. We define the following binary relation designated by $\varGamma$ on the set GN:
\begin{enumerate}[\textbullet]
\item $\cS\sim_{\varGamma}\cP$ iff $(\cS\sa\cP\wedge\cP\sa\cS)\in\varGamma$,
\end{enumerate}
We have:
\begin{enumerate}[\textbullet]
\item $\sim_{\varGamma}$ is an equivalence relation being a congruence with respect to all functors.
\end{enumerate}
\emph{The sketch of the proof}. By \eqref{Barbara} and \eqref{Ia}, $\sim_{\varGamma}$ is a equivalence relation and for all $\cS,\cP\in\mathrm{GN}$, we obtain: $\cS\sim_{\varGamma}\cP$ if and only if for every $\cM\in\mathrm{GN}$ we have: ($\cS\sa\cM\in\varGamma$ $\Leftrightarrow$ $\cP\sa\cM\in\varGamma$) and ($\cM\sa\cS\in\varGamma \Leftrightarrow \cM\sa\cP\in\varGamma$).

Moreover, by \eqref{Datisi} and \eqref{Ci}, for all $\cS,\cP\in\mathrm{GN}$, we obtain: if $\cS\sim_{\varGamma}\cP$, then for every $\cM\in\mathrm{GN}$ we have: ($\cS\si\cM\in\varGamma \Leftrightarrow \cP\si\cM\in\varGamma$) and ($\cM\si\cS\in\varGamma \Leftrightarrow \cM\si\cP\in\varGamma$).

Finally, by \eqref{Ish1}, \eqref{isSi}, \eqref{Datisi}, \eqref{Ci}, \eqref{cIi}, \eqref{9} and \eqref{8}, for all $\cS,\cP\in\mathrm{GN}$, we obtain: if $\cS\sim_{\varGamma}\cP$, then for every $\cM\in\mathrm{GN}$ we have: ($\cS\sis\cM\in\varGamma \Leftrightarrow\cP\sis\cM\in\varGamma$) and ($\cM\sis\cS\in\varGamma$ $\Leftrightarrow$ $\cM\sis\cP\in\varGamma$). \qed

Let $\|\cS\|$ be the equivalence class of \cS\ with respect to~$\sim_{\varGamma}$. By \eqref{Datisi}, \eqref{Ci} and \eqref{5}, for all $\cS,\cM\in\mathrm{GN}$, we get:
\begin{enumerate}[\textbullet]
\item if $\cS\sis\cS\in\varGamma$, $\cM\si\cM\in\varGamma$ and $\cM\sa\cS\in\varGamma$, then $\|\cM\|=\|\cS\|$.
\end{enumerate}

We use  $\mM_{\varGamma}=\langle\uU_{\varGamma}, \DD_{\varGamma}\rangle$, where:
\begin{enumerate}[\textbullet]
\item $\uU_{\varGamma}$ consists of all pairs $\{\cM,\cQ\}$ of name letters such that $\cM\si\cQ\in\varGamma$ and all equivalent classes $\|\cM\|$ and name letters $\cM$ such that $\cM\si\cM\in\varGamma$ ($\uU_{\varGamma}$ may be empty),
\item $\DD_{\varGamma}(\cS)\coloneqq
\begin{cases}
\{\|\cM\|\in\uU_{\varGamma} : \cM\sa\cS\in\varGamma\}=\{\|\cS\|\} &\text{if $\cS\sis\cS\in\varGamma$}\\
\{\{\cM,\cQ\}\in\uU_{\varGamma} :  (\cM\sa\cS\vee\cQ\sa\cS)\in\varGamma\} \cup\null &\\
\{\|\cM\|\in\uU_{\varGamma} : \cM\sa\cS\in\varGamma\}\cup \{\cM\in\uU_{\varGamma} : \cM\sa\cS\in\varGamma\} &\text{if $\cS\sis\cS\notin\varGamma$}
\end{cases}$
\end{enumerate}
In the case of $\cS\sis\cS\notin\varGamma$ we added a third set since it is possible that for each $\cM\in\mathrm{GN}$ we have $\cM\sa\cS\notin\varGamma$; and then $\|\cS\|=\{\cS\}$ and  $\{\|\cM\|\in\uU_{\varGamma} : \cM\sa\cS\in\varGamma\}=\{\|\cS\|\}=\{\{\cS\}\}=\{\{\cM,\cQ\}\in\uU_{\varGamma} :  (\cM\sa\cS\vee\cQ\sa\cS)\in\varGamma\}$.

For every name letter $\cS$, we obtain:
\begin{enumerate}[(i)]
\item\label{cbf} if $\cS\sis\cS\notin\varGamma$ and $\cS\si\cS\in\varGamma$, then $\{\cS,\{\cS\},\|\cS\|\}\subseteq\DD_{\varGamma}(\cS)$;
\item\label{abf} $\cS\si\cS\in\varGamma$ iff   $\|\cS\|\in\DD_{\varGamma}(\cS)$ iff   $\DD_{\varGamma}(\cS)\neq\emptyset$.

\item\label{dbf} $\cS\sis\cS\in\varGamma$ iff  $\DD_{\varGamma}(\cS)=\{\|\cS\|\}$ iff $\DD_{\varGamma}(\cS)$ is a singleton.
\end{enumerate}
Indeed, for \eqref{cbf}: By \eqref{Ia}. For \eqref{abf}: By \eqref{cbf} and using \eqref{cIi}, \eqref{Ci} and \eqref{Daraptiplus}. For \eqref{dbf}: If $\DD_{\varGamma}(\cS)$ is a singleton, then $\cS\sis\cS\in\varGamma$, by \eqref{cbf} and~\eqref{abf}.

By induction, for any formula $\alpha$, we obtain:
\begin{enumerate}[(C\Shis)]
\item $\alpha$ is true in $\mM_{\varGamma}$ iff $\alpha\in\varGamma$.
\end{enumerate}
\emph{The sketch of the proof}. For ${\circ}\in\{\sa,\si,\se,\so,\sis\}$ and all $\cS,\cP\in\mathrm{GM}$, we have:
\begin{enumerate}[\textbullet]
\item $\cS\circ\cP$ is true in $\mM_{\varGamma}$ iff $\cS\circ\cP\in\varGamma$.
\end{enumerate}

As for \Sh, we obtain: if $\cS\sa \cP$ is true in $\mM_{\varGamma}$, then   $\cS\sa\cP\in\varGamma$. For the proof of the converse implication we assume that $\cS\sa\cP\in\varGamma$. We will consider three cases. The first one: $\cS\si\cS\notin\varGamma$. Then, by \eqref{abf},  $\emptyset=\DD_{\varGamma}(\cS) \subseteq\DD_{\varGamma}(\cP)$. The second case: $\cP\sis\cP\notin\varGamma$. If $\cM\in\DD_{\varGamma}(\cS)$ (resp.\ $\|\cM\|\in\DD_{\varGamma}(\cS)$,  $\{\cM,\cQ\}\in\DD_{\varGamma}(\cS)$), then $\cM\in\DD_{\varGamma}(\cP)$ (resp.\ $\|\cM\|\in\DD_{\varGamma}(\cS)$,  $\{\cM,\cQ\}\in\DD_{\varGamma}(\cP)$), by \eqref{Barbara}. The third case: $\cS\si\cS\in\varGamma$ and $\cP\sis\cP\in\varGamma$. Then, by \eqref{6}, $\cS\sis\cS\in\varGamma$. Moreover, by \eqref{Datisi} and \eqref{Ci}, $\cS\si\cP\in\varGamma$ and $\cP\si\cS\in\varGamma$. Hence and from \eqref{5}, also $\cP\sa\cS\in\varGamma$. So, $\|\cS\|=\|\cP\|$ and $\DD_{\varGamma}(\cS)=\DD_{\varGamma}(\cP)$.  So, in all three cases, $\DD_{\varGamma}(\cS)\subseteq\DD_{\varGamma}(\cP)$.

Let $\cS\si\cP$ be true in $\mM_{\varGamma}$, i.e.\ some $\cM$ or $\|\cM\|$ or $\{\cM,\cQ\}$ belongs to $\DD_{\varGamma}(\cS)$ and $\DD_{\varGamma}(\cP)$. In the first and second cases, $(\cM\si\cM\wedge\cM\sa\cS\wedge\cM\sa\cP)\in\varGamma$. So $\cS\si\cP\in\varGamma$, by \eqref{Daraptiplus}. In the third case, the same as for \Sh. For the proof of the converse implication we assume that $\cS\si\cP\in\varGamma$. Then, by \eqref{cIi} and \eqref{Ci}, also $\cS\si\cS\in\varGamma$, $\cP\si\cS\in\varGamma$ and $\cP\si\cP\in\varGamma$.  We will consider three cases. The first one: $\cS\sis\cS\in\varGamma$. Then, by \eqref{5},  $\cS\sa\cP\in\varGamma$. Hence $\|\cS\|\in\DD_{\varGamma}(\cP)$ and $\DD_{\varGamma}(\cS)=\{\|\cS\|\}\subseteq\DD_{\varGamma}(\cP)$. The second case: $\cS\sis\cS\notin\varGamma$ and $\cP\sis\cP\in\varGamma$. Then, by \eqref{5},  $\cP\sa\cS\in\varGamma$. Hence $\DD_{\varGamma}(\cP)=\{\|\cP\|\}\subseteq\DD_{\varGamma}(\cS)$. The third case: $\cS\sis\cS\notin\varGamma$ and $\cP\sis\cP\notin\varGamma$. It is the same as for \Sh. So, in all three cases, $\DD_{\varGamma}(\cS)\cap\DD_{\varGamma}(\cP)\neq\emptyset$.

Now, let $\cS\sis\cP$ be true in $\mM_{\varGamma}$, i.e., $\DD_{\varGamma}(\cS)$ be a singleton whose only element belongs to~$\DD_{\varGamma}(\cP)$. Then, $\DD_{\varGamma}(\cS\subseteq\DD_{\varGamma}(\cP)$, and so   $\cS\sa\cP\in\varGamma$. Moreover, $\cS\sis\cS\in\varGamma$, by \eqref{dbf}. Hence $\cS\sis\cP\in\varGamma$, by \eqref{4}. For the proof of the converse implication, let $\cS\sis\cP\in\varGamma$. Then $\cS\sis\cS\in\varGamma$, by \eqref{Ish1}; and so $\DD_{\varGamma}(\cS)=\{\|\cS\|\}$. Moreover, $\cS\sa\cP\in\varGamma$ and $\cS\si\cS\in\varGamma$, by \eqref{isSa} and \eqref{isSi}. Hence $\|\cS\|\in\DD_{\varGamma}(\cP)$. So, $\cS\sis\cP$ is true in~$\mM_{\varGamma}$.

The rest of the proof is similar to that for \Sh. \qed

Therefore, by \eqref{CShis}, using the properties of maximal consistent sets,  we get:
\begin{enumerate}[\textbullet]
\item \emph{A formula is a thesis of\/ \textup{\Shis} if and only if it is an\/ \sa\si\se\so\sis-tautology}.
\end{enumerate}
\begin{Remark}
For \Shis, such as for \Sh\, we can use only models with a non-empty universe. Indeed,  it can also be shown by taking canonical models with non-empty $\uU_{\varGamma}$:
\begin{enumerate}[\textbullet]
\item $\uU_{\varGamma}$ consists of all pairs of name letters, all equivalent classes of name letters and all name letters,
\item $\DD_{\varGamma}(\cS)\coloneqq
\begin{cases}
\{\|\cM\| : \cM\si\cM\wedge\cM\sa\cS\in\varGamma\}=\{\|\cS\|\} &\text{if $\cS\sis\cS\in\varGamma$}\\
\{\{\cM,\cQ\} :  \cM\si\cQ\wedge (\cM\sa\cS\vee\cQ\sa\cS)\in\varGamma\} \cup\null &\\
\{\|\cM\| : \cM\si\cM\wedge\cM\sa\cS\in\varGamma\}\cup &\\
\{\cM : \cM\si\cM\wedge\cM\sa\cS\in\varGamma\} &\text{if $\cS\sis\cS\notin\varGamma$}
\end{cases}$
\end{enumerate}
As for the previous model, we show that condition \eqref{CShis} holds.\qed
\end{Remark}
\begin{Remark}
Once again we see that defining $\uU_{\varGamma}$ for $\Sh$ and $\Lcal$ as for $\Shis$ and defining $\DD_{\varGamma}$ for these first calculi using the second of the alternative conditions used for $\Shis$ gives $\DD_{\varGamma}$ with the set of values without singletons. In this condition, the set of equivalent classes can be omitted. Just use the sum of two sets. Hence we get that for \Sh\ and \Lcal\ we can use non-monoreferential and polyreferential semantics, respectively. \qed
\end{Remark}

\section{Further extensions of calculi of names\label{sec7}}

In this section, we briefly present other possible extensions of the systems considered in the paper by adding new kinds of singular sentences and identities.

\subsection{Calculi of names plus traditional singular sentences and identities\label{subsec7.1}}

We can extend all of the calculi of names considered earlier to include traditional singular sentences and identities, which we discussed in Remark~\ref{Rem5.1}.  We remember that these sentences have the form: ‘$a$ is a $P$’, ‘$a$ is not a $P$’, and ‘$a$ is identical to $b$’, respectively. Symbolically, we will write them as ‘$a\sist P$’, ‘$a\bar{\sist} P$’ and ‘$a\idsf b$’. We also remember that only names with exactly one referent can be inserted for ‘$a$’ and ‘$b$’ (for ‘$P$’, we can use any general name).

To the set GN of general name letters, we add the countably infinite set SN of singular name letters (for which we use ‘$a$’, `$b$' and ‘$c$’ with or without indices). We build the new set of sentence formulas in the standard way.

Now, we have to use models with an additional denotation function for the singular letters. So these models will have the form $\la\uU,\DD,\dd\ra$, where $\uU$ is a non-empty set and $\dd$ is a function which assigns to any singular name letter an element of $\uU$. Using the natural interpretation of the functors ‘$\sist$’, ‘$\bar{\sist}$’ and ‘$\idsf$’, we extend the notions of \emph{being a true formula} in a model $\mM=\la\uU,\DD,\dd\ra$. For all $\ea,\eb\in\mathrm{SN}$ and $\cP\in\mathrm{GN}$, we assume:
\begin{enumerate}[\textbullet]
\item $\ea\sist\cP$ is true in $\mM$ iff  $\dd(\ea)\in\DD(\cP)$;
\item $\ea\bar{\sist}\cP$ is true in $\mM$ iff  $\dd(\ea)\notin\DD(\cP)$;
	\item $\ea\idsf\eb$ is true in $\mM$ iff $\dd(\ea)=\dd(\eb)$.
\end{enumerate}

Now, we extend all of the calculi of names considered earlier by adding the following tautologies as their additional axioms:\vspace{-3pt}
\begin{gather*}	
(a\sist M\wedge M\sa P)\to a\sist P\\
(a\sist M\wedge a\sist P)\to M\si P\\
a\bar{\sist} P\leftrightarrow \neg\, a\sist P\\
a\idsf a\\
a\idsf b \to b\idsf a\\
(a\idsf c \wedge  c\idsf b)\to a\idsf b\\
(a\idsf b \wedge a\sist P)\to b\sist P
\end{gather*}
The last four make ‘$\idsf$’ an equivalence relation being a congruence with respect to ‘$\sist$’.

Using the above axioms and suitable definitions, we get:
\begin{gather*}	
(a\sist M\wedge M\se P)\to a\bar{\sist} P\\
(a\sist M\wedge a\bar{\sist} P)\to M\so P
\end{gather*}
\indent In conclusion, by appropriately applying Henkin's method, we can prove the completeness of the extended versions of the calculi of names studied.

\subsection{Calculi of names plus Czeżowski’s singular sentences and identities\label{subsec7.2}}

Tadeusz \citet{CzM} analysed singular sentences with a subject of the form ‘this $S$’, where ‘$S$’ is to be replaced by a   non-empty general name (cf.\ Remark~\ref{Rem5.1}). He assumed that: “The name, ‘This~$S$’, in the subject of a singular proposition I regard to be a proper name denoting a given individual from the extension of the S term’’ \citep[(p.~392]{CzM}. Therefore, for Czeżowski ‘This $S$ is an $S$’ is a tautology; in symbolic notation:
\[\label{It}\tag{I\st}
\st S\sist S
\]
Tautologies are also all seven formulas that we obtain from the axioms from the previous point by substitution: $a/\st S$, $b/\st M$ and $c/\st P$.

We extend models used for Łukasiewicz’s calculus \Lcal, adding a choice function $\dc$, which for any general name letter “indicates” one of its referents; i.e., $\dc(\cS)\in\DD(\cS)$ for every name letter $\cS$. Using the natural interpretation of the functors ‘$\sist$’, ‘$\bar{\sist}$’ and ‘$\idsf$’, we extend the notions of \emph{being a true formula} in a model $\mM=\la  \uU,\DD,\dc\ra$. For all $\cS,\cP\in\mathrm{GN}$, we assume:
\begin{enumerate}[\textbullet]
\item $\st\cS\sist\cP$ is true in $\mM$ iff $\dc(\cS)\in\DD(\cP)$;
\item $\st\cS\bar{\sist}\cP$ is true in $\mM$ iff $\dc(\cS)\notin\DD(\cP)$;
\item $\st\cS\idsf\st\cP$ is true in $\mM$ iff $\dc(\cS)=\dc(\cP)$.
\end{enumerate}

Using Henkin's method, we can prove that by adding eight new axioms to \Lcal, i.e.\ \eqref{It} and seven tautologies which we obtain from the axioms from \eqref{subsec7.1} by substitution, we obtain a complete calculus with respect to the above semantics.

One can ask the following question:
\begin{enumerate}[\textbullet]
\item What happens if we reject Czeżowski’s assumption that the object chosen as $S$ is an~$S$?
\end{enumerate}
We may indicate a given object as one of the $S$s, but it is not. When we allow this, \eqref{It} ceases to be a tautology. In the described situation, however, the following additional problem arises:
\begin{enumerate}[\textbullet]
\item Are ‘This $S$ is a $P$’ and ‘This $S$ is not a $P$’ without truth values or only false?
\end{enumerate}
Depending on the answer, different logical systems (two-valued and three-valued) can be created.

\renewcommand{\bibfont}{\small}
\setlength{\bibsep}{3pt}


\begin{thebibliography}{99}
\bibitem[Clark(1980)]{Clark} Clark, M., 1980, \emph{The Place of Syllogistic in Logical Theory}, Nottingham University Pres: Nottingham.


\bibitem[Corcoran(1972a)]{Cor72a} Corcoran, J., 1972a, ``Aristotle's natural deduction system'', \emph{Journal of Symbolic Logic}, 37: 437 (Abstract).

\bibitem[Corcoran(1972b)]{Cor72b} Corcoran, J., 1972b, ``Completeness of an ancient logic'', \emph{Journal of Symbolic Logic}, 37: 696--702. \url{https://doi.org/10.2307/2272415}


\bibitem[Corcoran(1974)]{Cor74} Corcoran, J., 1974,, ``Aristotle's natural deduction system'', pages 85--131 in J.~Corcoran (ed.), \emph{Ancient Logic and Its Modern Interpretations}, Reidel: Dordrecht. \url{https://doi.org/10.1007/978-94-010-2130-2_6}




\bibitem[Czeżowski(1955)]{CzM} Czeżowski, T., 1955, ``On certain peculiarities of singular propositions'', \emph{Mind}, 255: 392--395. \url{https://doi.org/10.1093/mind/LXIV.255.392}

\bibitem[Ishimoto(1977)]{Ish} Ishimoto, A., 1977, ``A propositional fragment of Leśniewski's elementary ontology'', \emph{Studia Logica},  36: 285--299. \url{https://doi.org/10.1007/BF02120666}

\bibitem[Kotarbiński(1966)]{Kot} Kotarbiński, T., 1966, \emph{Gnosiology.\ The Scientific Approach to the Theory of Knowledge}; Pergamon Press and Ossolineum: Oxford and Wrocław. English version of the second Polish revised edition: \emph{Elementy teorii poznania, logiki formalnej i metodologii nauk}, Ossolineum: Wrocław, Poland, 1961. The first Polish edition: Zakład Narodowy im.\ Ossolińskich: Lwów, Poland, 1929.

\bibitem[Lejewski(1984)]{Lej} Lejewski, C., 1984, ``On Leśniewski's ontology'', pages 1230--148 in J.\ts T.\ts J. Srzednicki and V.\ts F. Rickey (eds.), \emph{Leśniewski's Systems.\ Ontology and Mereology}, Martinus Nijhoff Publishers (Kluwer) and Ossolineum: The Hague and Wrocław. Reprinted from: \emph{Ratio} \textbf{1958}, \textit{1}, 150--176.


\bibitem[Łukasiewicz(1934)]{Luk34} Łukasiewicz, J., 1934, ``Znaczenie analizy logicznej dla poznania'' (The importance of logical analysis for cognition), \emph{Przegląd Filozoficzny}, 37: 369--377. It is the translation of the German lecture entitled ``Bedeutung der logischen Analyse für die Erkenntnis'', delivered at the ``VIII Congrès International de Philosophie'', Prague, September 2--7, 1934. The original was published in \emph{Actes du VIII Congrès International de Philosophie}, Prague, Czechia, 1936, pp.~75--84.

\bibitem[Łukasiewicz(1939)]{Luk39} Łukasiewicz, J., 1939, ``O sylogistyce Arystotelesa'' (On Aristotle's syllogistic), \emph{Sprawozdania PAU}, 44: 220--227. Reprint in: Łukasiewicz, J. \emph{Z~zagadnień logiki i filozofii} (From the Issues of Logic and Philosophy), PWN: Warszawa, Poland, 1961, pp.~220--227. It summarises the lecture presented by Łukasiewicz on June 9, 1939, at a Polish Academy of Arts and Sciences (PAU) session. He presented his fully completed monograph on Aristotle’s syllogistic in this lecture. This monograph was destroyed at the beginning of World War~II. After the war, in Dublin, Łukasiewicz recreated his results, publishing 1951 the monograph \citep{Luk57}.

\bibitem[Łukasiewicz(1957)]{Luk57} Łukasiewicz, J., 1957, \emph{Aristotle's Syllogistic from the Standpoint of Modern Formal Logic}, 2nd ed.,
    Clarendon Press: Oxford. The first edition:  Oxford Univ.\ Press.: Oxford, 1951.


\bibitem[Łukasiewicz(1963)]{Luk63} Łukasiewicz, J., 1963, \emph{Elements of Mathematical Logic}, The Macmillan Co.: New York. English transl.\ by O.~Wojtasiewicz of the second Polish edition: \emph{Elementy logiki matematycznej}, Polish Scientific Publ.: Warsaw, 1958. The first Polish edition:  ed.\ by M.~Presburger, Publications of Students of Mathematics and Physics in Warsaw Univesity, vol.~18: Warsaw,  1929.


\bibitem[Morawiec(1961)]{M} Morawiec, A., 1961, ``Podstawy logiki nazw'', \emph{Studia Logica}, 12: 145--161. Summary: ``Foundations of the theory of names'', \emph{Studia Logica}, 12: 167--170. \url{https://doi.org/10.1007/BF02126831}



\bibitem[Moss(2008)]{Moss} Moss, L.\ts S., 2008, ``Completeness theorems for syllogistic fragments'', pages 143--173 in F.~Hamm and S.~Kepser (eds.), \emph{Logics for Linguistic Structures}, Mouton de Gruyter: Berlin and New York. \url{https://doi.org/10.1515/9783110211788.143}


\bibitem[Pietruszczak(1987)]{ja87} Pietruszczak, A., 1987, ``O logice tradycyjnej i rachunku nazw dopuszczającym podstawienia nazw pustych'' (On traditional logic and calculus of names allowing substitutions of empty names), \emph{Ruch Filozoficzny}, 44: 158--166.

\bibitem[Pietruszczak(1991a)]{ja91a} Pietruszczak, A., 1991a, ``Standardowe rachunki nazw z funktorem Leśniewskiego'' (Standard calculi of names with Leśniewski's functor). \emph{AUNC Logika}, 1: 5--29.


\bibitem[Pietruszczak(1991b)]{ja91b} Pietruszczak, A., 1991b, ``O pewnym ujęciu logiki tradycyjnej'' (On a certain approach to traditional logic), \emph{AUNC Logika}, 1: 31--41.

\bibitem[Pietruszczak(1991c)]{ja91c} Pietruszczak, A., 1991c, \emph{Bezkwantyfikatorowy rachunek nazw.\ Systemy i ich metateoria} (Quantifier-free calculus of names.\ Systems and their metatheory), Wydawnictwo Adam Mar\-sza\-łek: Toruń. It is a book edition of a part of the doctoral dissertation, Nicolaus Copernicus University in Toruń, 1990.


\bibitem[Pietruszczak(1992)]{ja92} Pietruszczak, A., 1992, `` Stała Leśniewskiego w teoriach sylogistycznych.\ Semantyczne badania pewnych kwantyfikatorowych rachunków nazw'' (Leśsniewski's copula in syllogistic theories.\ Semantic studies of certain quantifier calculus of names), \emph{AUNC Logika}, 3: 45--76.


\bibitem[Pratt-Hartmann and Moss(2009)]{P-HiM} Pratt-Hartmann, I., and L.\ts S.~Moss, 2009, ``Logics for the relational syllogistic'', \emph{The Review of Symbolic Logic}, 2(4): 647--683. \url{https://doi.org/10.1017/S1755020309990086}

\bibitem[Quine(1950)]{Q50} Quine, W.\ts V.\ts O., 1950, \emph{Methods of Logic}; Henry Holt and Co.: New York.


\bibitem[Quine(1953)]{Q53} Quine, W.\ts V., 1953, ``Logic and reification of universals'', pages 102--129 in \emph{From a Logical Point of View}, Harvard University Press.

\bibitem[Rasiowa and Sikorski(1968)]{RiS} Rasiowa, H. and R.~Sikorski, 1968, \emph{The Mathematics of Metamatematics}, 2nd ed., PWN: Warszawa.


\bibitem[Shepherdson(1956)]{Sh} Shepherdson, J.\ts C., 1956, ``On the interpretation of Aristotelian syllogistic'', \emph{The Journal of Symbolic Logic}, 21: 137--147. \url{https://doi.org/10.2307/2268752}

\bibitem[Słupecki(1946)]{Sl} Słupecki, J., 1946, ``Uwagi o sylogistyce Arystotelesa'' (Notes on Aristotle's syllogistic), \emph{Annales UMCS}, sectio~F, 1: 187--191.

\bibitem[Smiley(1973)]{Smil} Smiley, T.\ts J., 1973, ``What is a syllogism?'', \emph{Journal of Philosophical Logic}, 2: 136--154. \url{https://doi.org/10.1007/BF02115614}

\bibitem[Smith(1989)]{Sm} Smith, R., 1989, ``Introduction'', pages XIII--XXXI in Aristotle, \emph{Prior Analytics}, translated, with introduction, notes, and commentary, by R.~Smith, Hackett Publishing Company: Indianapolis, Cambridge.

\bibitem[Takano(1985)]{T} Takano, M., 1985, ``A semantical investigation into Leśniewski's axiom of his ontology'', \emph{Studia Logica}, 44: 71--77.  \url{https://doi.org/10.1007/BF00370810}

\end{thebibliography}
\end{document}